\documentclass[12pt]{amsart}

\usepackage{amsfonts, amsmath, amssymb}
\usepackage{graphicx, graphics}

\newtheorem{prop}{Proposition}[section]
\newtheorem{theo}[prop]{Theorem}
\newtheorem{coro}[prop]{Corollary}

\newtheorem{conj}[prop]{Conjecture}

\theoremstyle{definition}
\newtheorem{dfn}[prop]{Definition}
\newtheorem{exam}[prop]{Example}
\theoremstyle{remark}
\newtheorem{rem}[prop]{Remark}

\def\<{\langle}
\def\>{\rangle}

\def\Si{{\Sigma}}

\def\O{{\mathcal O}}
\def\P{{\mathbb P}}
\def\R{{\mathbb R}}
\def\T{{\mathbb T}}
\def\C{{\mathbb C}}
\def\Z{{\mathbb Z}}
\def\Q{{\mathbb Q}}

\begin{document}

\bigskip

\title{Toric Residues and Mirror Symmetry}

\author{Victor~V.~Batyrev}
\address{Mathematisches Institut, Universit\"at T\"ubingen,
Auf der Morgenstelle 10,\newline T\"ubingen D-72076, Germany} 
\email{victor.batyrev@uni-tuebingen.de} 

\author{Evgeny~N.~Materov}
\address{Mathematisches Institut, Universit\"at T\"ubingen,
Auf der Morgenstelle 10,\newline T\"ubingen D-72076, Germany} 
\email{evgeny.materov@uni-tuebingen.de}


\dedicatory{To Yuri Ivanovich Manin on his 65-th birthday}


\begin{abstract}
We develop some ideas of Morrison and Plesser and formulate a precise
mathematical conjecture which has close relations to toric mirror
symmetry. Our conjecture, we call it Toric Residue Mirror Conjecture, claims
that the generating functions of intersection numbers of divisors on a special
sequence of simplicial toric varieties are power series expansions of some
rational functions obtained as toric residues. We expect that this conjecture
holds true for all Gorenstein toric Fano varieties associated with reflexive
polytopes and give some evidences for that. The proposed conjecture suggests a
simple method for computing Yukawa couplings for toric mirror Calabi-Yau
hypersurfaces without solving systems of differential equations. We make
several explicit computations for Calabi-Yau hypersurfaces in weighted
projective spaces and in products of projective spaces. 
\end{abstract}

\maketitle

\tableofcontents

\thispagestyle{empty}

\newpage

\section{Introduction} 

The mirror symmetry attracts interest of mathematicians because it allows to
identify mathematical objects of a very different nature: generating functions
for Gromov-Witten invariants of rational curves on Calabi-Yau manifolds $X$
and power series expansions of special functions on periods of  the mirror
family of Calabi-Yau manifolds $X^*$. Many examples of this identification can
be computed  explicitly for Calabi-Yau hypersurfaces in toric Fano
varieties. The general toric mirror construction \cite{Batyrev2} suggests a
duality between Calabi-Yau varieties with Gorenstein canonical
singularities. These singularities in general can not be resolved without
changing the canonical class. Therefore, a verification of the toric mirror
symmetry in full generality requires orbifold versions  of quantum cohomology
and Gromov-Witten invariants for singular varieties (cf. \cite{CR1,CR2}). The
Mirror Theorem states that power series obtained by these two very different
methods are actually the same. It is rather nontrivial to prove the  Mirror
Theorem even in the very special case of Calabi-Yau quintic $3$-folds
\cite{G2,G4}. 
\smallskip

In this paper we propose a more elementary mirror symmetry test which is
expected to hold for {\it all} families of Calabi-Yau hypersurfaces in
Gorenstein toric Fano varieties associated with dual pairs of reflexive
polytopes. We formulate this test as a mathematical conjecture and call it
{\it Toric Residue Mirror Conjecture}. The idea of this conjecture is due to
Morrison and Plesser \cite{MP} who have checked it for some toric Calabi-Yau
hypersurfaces (including the case of Calabi-Yau quintic
$3$-folds). Unfortunately, Morrison and Plesser didn't formulate their ideas
in the form of a precise mathematical conjecture. The most important
contribution of Morrison and Plesser \cite{MP} is the construction of a
sequence of toric manifolds $\{\P_{\beta}\}$ associated with pairs 
$(\P,\beta)$,  where $\P$ is a  smooth projective $d$-dimensional toric
manifold $\P$ and  $\beta$ is an element in the monoid of integral points in
the Mori cone $K_{\rm eff}(\P)$ of $\P$. Using ideas of Morrison and Plesser,
we define a special  cohomology class 
$\Phi_{\beta} \in H^{2(\dim \P_\beta -d)}(\P_{\beta}, \Q)$ which plays the
role of a ``virtual fundamental class'' in our mirror symmetry test. The
manifold $\P_{\beta}$ and the cohomology class $\Phi_{\beta}$ will be called
{\it Morrison-Plesser moduli space} and {\it Morrison-Plesser class}
respectively. For our Toric Residue Mirror Conjecture, we will need a slight
generalization of the Morrison-Plesser construction for projective simplicial
toric varieties $\P$. 
\smallskip

Let us make some remarks concerning the relation of our conjecture to the
traditional point of view on the  mirror symmetry for Calabi-Yau manifolds
(see, e.g., \cite{CK}). First of all, it is important to emphasize that Toric
Residue Mirror Conjecture can not be obtained as a consequence of the mirror
theorems due to Givental and Lian-Liu-Yau \cite{G4,LLY2,LLY3}. On the other
hand, it seems that all known methods for proving Mirror Theorem for toric
Calabi-Yau hypersurfaces and complete intersections use some versions of the
Morrison-Plesser moduli spaces $\P_{\beta}$ as naive toric approximations of
the Kontsevich moduli spaces of stable maps
\cite{G1,G2,G3,G4,G5,LLY1,LLY2,LLY3}. 

This makes impression that our conjecture could provide a helpful intermediate
step for future formulations and verifications of the toric mirror symmetry
predictions in case of Calabi-Yau varieties with orbifold singularities. We
remark that the  Toric Residue Mirror Conjecture also suggests a simple method
for computing Yukawa $(d-1)$-point functions for $(d-1)$-dimensional toric
Calabi-Yau hypersurfaces  without  using systems of differential equations. 
\medskip

The main advantage of the Toric Residue Mirror Conjecture is its
simplicity. One does not need to know much (e.g., Hodge theory,
Gromov-Witten invariants and quantum cohomology) in order to understand its 
statement. Our conjecture uses only the intersection theory on simplicial
toric varieties $\P_{\beta}$ and the notion of toric residues which are
standard operations in computational commutative algebra. Unfortunately, we
were not able to prove our conjecture in general. In this paper, we check it
for some  classes of reflexive polytopes using direct computations. We hope
that a conceptual proof of the Toric Residue Mirror Conjecture for arbitrary
reflexive polytopes would be an important  contribution  to understanding of
the  mirror symmetry phenomenon. 
\bigskip 

The paper is organized as follows. In Section 2, we give a review of toric
residues and their properties. In Section 3, we discuss Morrison-Plesser
moduli spaces $\P_{\beta}$ associated with lattice points $\beta$ in the Mori
cone $K_{\rm eff}(\P)$ of a simplicial toric variety $\P$. In Section 4, we
formulate the Toric Residue Mirror Conjecture and give some evidences for
it. In Sections 5 and 6, we illustrate our conjecture for some simplest
interesting cases. In Sections 7 and 8, we prove the Toric Residue Mirror
Conjecture for Calabi-Yau hypersurfaces in weighted projective spaces and in
product of projective spaces. Section 9 explains applications of our
conjecture to the toric mirror symmetry and to the computations of Yukawa
couplings for  toric Calabi-Yau hypersurfaces. Some examples of these
computations are given in Sections 10 and 11. 
\medskip

We remark that the Toric Residue Mirror Conjecture can be  formulated in the
same generality also for toric Calabi-Yau complete intersections associated
with nef-partitions of reflexive polytopes. The details of this generalization
will be published in the forthcoming paper \cite{BM}.  

\bigskip

{\it Acknowledgments.} We would like to thank A.~Mavlyutov for suggesting us
to use his formula for computation of Yukawa couplings and for pointing out us
on the correction in constants in this formula. We thank E.~Cattani,
A.~Dickenstein and M.~Passare for their remarks and corrections concerning
preprint version of this work. The authors were supported by DFG,
Forschungsschwerpunkt ``Globale Methoden in der komplexen
Geometrie''. E.~Materov was partially supported by RFBR Grant 00-15-96140. 

\bigskip

\section{Toric residues}
\label{Sect_Residues}

The notion of toric residue was introduced by D.~Cox using homogeneous
coordinates on toric varieties \cite{Cox2}. In this paper, we prefer 
another approach which uses
the affine coordinates $t_1, \ldots, t_d$ on the algebraic torus 
$\T = {\rm Spec}\, \C [t_1^{\pm 1}, \ldots, t_d^{\pm 1}]$. 
\medskip

Denote by $M$ the group of characters of $\T$ which is isomorphic to
$\Z^d$. Let $\Delta$ be a $d$-dimensional convex polytope in 
$M_{\R} = M \otimes \R$ with  vertices in $M$. 

\begin{dfn}
We denote by 
\[
  S_{\Delta} = \bigoplus_{k = 0}^{\infty} S_{\Delta}^k 
\]
the graded subring in $\C[t_0, t_1^{\pm 1}, \ldots, t_d^{\pm 1}]$ whose
$k$-th homogeneous component $S_{\Delta}^k$ is spanned as $\C$-vector space by
all monomials $t_0^k t_1^{m_1} \cdots t_d^{m_d}$ such that the lattice point 
$(m_1, \ldots, m_d)$ is contained in $k\Delta$. 
\end{dfn}

\begin{dfn}
The variety $\P_{\Delta}:= {\rm Proj}\, S_{\Delta}$ is called 
 {\it the projective toric variety}  
{\it associated with the polytope} $\Delta$. We denote by 
$\O_{\P_\Delta}(1)$ the ample invertible sheaf
on $\P_\Delta$ corresponding to the graded $S_\Delta$-module $S_\Delta(-1)$.
\end{dfn}

\begin{dfn}
Denote by $I_{\Delta} = \bigoplus_{k =0}^{\infty} I_{\Delta}^k$ the graded
ideal in $S_{\Delta}$ whose $k$-th homogeneous component $I_{\Delta}^k$ is
spanned over $\C$ by all monomials $t_0^k t_1^{m_1} \cdots t_d^{m_d}$ such
that $(m_1, \ldots, m_d)$ is contained in the interior of $k\Delta$. 
\end{dfn}

It is known that $S_\Delta$ is a Cohen-Macaulay ring and $I_\Delta$ is its
dualizing module \cite{Danilov}.
Let $g_0(t),g_1(t),\ldots,g_d(t)$ be generic Laurent polynomials in 
$\C [t_1^{\pm 1}, \ldots, t_d^{\pm 1}]$ with supports in $\Delta$ such that 
\[
  G = (G_0, G_1,   \ldots, G_d) := 
      (t_0 g_0(t), t_0 g_1(t), \ldots, t_0 g_d(t))\subset S^1_{\Delta}
\]
is a regular sequence in $S_{\Delta}$. We can consider the polynomials 
$G_0, G_1,\ldots, G_d$ also as global sections of 
$\O_{\P_\Delta}(1)$ on $\P_\Delta$ having  no common zeros. 

Since $G$ is also a regular sequence  for the dualizing module $I_\Delta$, 
we obtain two finite-dimensional graded $\C$-vector spaces    
\[
  S_{G} := S_{\Delta}/\< G_0, G_1,\ldots,G_d \>S_{\Delta}, \quad
  I_{G} := I_{\Delta}/\< G_0, G_1,\ldots,G_d \>I_{\Delta}, 
\]
where  $I_{G} $ is a dualizing module of 
the Artinian graded ring $S_{G}$ (see 
\cite[Proposition~9.4]{Batyrev1}). One has   a duality 
\[
  S^k_G \times I_G^{d + 1 - k} \to I_G^{d + 1}, \quad k = 0, \ldots, d,
\]
where  $I^{d + 1}_G$ is a $1$-dimensional $\C$-vector space. The Jacobian 
\begin{equation}
\label{toric_Jac}
  J_G:= \det 
  \left(t_i \frac{\partial G_j}
  {\partial t_i} \right)_{0 \leq i, j \leq d} = 
  t_0^{d + 1} \det
  \begin{pmatrix} 
  g_0                                   &  
  g_1                                   &  
  \cdots                                & 
  g_d                                   \\
  t_1 \frac{\partial g_0}{\partial t_1} & 
  t_1 \frac{\partial g_1}{\partial t_1} & 
  \cdots                                & 
  t_1 \frac{\partial g_d}{\partial t_1} \\
  \vdots &  
  \vdots & 
  \ddots & 
  \vdots  \\
  t_d \frac{\partial g_0}{\partial t_d} & 
  t_d \frac{\partial g_1}{\partial t_d} & 
  \cdots & 
  t_d \frac{\partial g_d}{\partial t_d}
  \end{pmatrix} 
\end{equation}
is an element of $S^{d + 1}_{\Delta}$. One can easily check that 
$J_G$
is contained in  $I^{d + 1}_{\Delta} \subset S^{d + 1}_{\Delta}$ 
(see \cite[Proposition 1.2]{CDS}) 
and  the class of $J_G$ in $I_G$ spans $I^{d + 1}_G$.

We define the toric residue in the following algebraic way: 

\begin{dfn} 
Let ${\rm Vol}(\Delta)$ be the {\it normalized volume} of
$\Delta$, i.e., $d!\cdot$(volume of $\Delta)$. The $\C$-linear map 
\begin{equation}
\label{Tor_res_map}
  {\rm Res}_G : I^{d + 1}_{\Delta} \rightarrow \C, 
\end{equation}
vanishing on the $\C$-subspace $\< G_0, G_1, \ldots, G_d\>I^{d}_{\Delta} 
\subset I^{d + 1}_{\Delta}$ and sending $J_{G}$
to ${\rm Vol}(\Delta)$ is called the {\bf toric residue}. 
This map establishes a canonical isomorphism 
\[
  I_G^{d + 1} = 
  I_\Delta^{d + 1}/\< G_0, G_1, \ldots, G_d\>I^{d}_{\Delta}\cong \C.
\]
\label{Tor_res}
\end{dfn} 

\begin{rem} This algebraic definition of the toric residue works as well for
any algebraically closed field $k$ of characteristic $0$. 
 We compare it with the 
definition given by D.~Cox in Section~\ref{Sect_Yukawa_CY}. 
\end{rem}

There exists a more geometric point of view on the toric residue. 
For this, one  remarks that the coherent sheaf $\widetilde{I_{\Delta}}$
associated with 
the graded $S_{\Delta}$-module $I_{\Delta}$ is exactly the sheaf 
$\Omega^d_{\P_{\Delta}}$ of Zariski differential $d$-forms (or dualizing
sheaf) on $\P_{\Delta}$ \cite{Danilov}. Comparing the pairing 
\[
  H^d(\P_{\Delta},\Omega^d_{\P_{\Delta}})\times 
  H^0(\P_{\Delta},\O_{\P_{\Delta}}) \rightarrow
  H^d(\P_{\Delta},\Omega^d_{\P_{\Delta}})
\]
with the pairing 
\[
  I_G^{d + 1}\times S_G^0 \rightarrow I_G^{d + 1},
\]
we can canonically identify the $1$-dimensional $\C$-space $I_G^{d + 1}$ with 
$H^d(\P_{\Delta},\Omega^d_{\P_{\Delta}})$. The regular sequence $G$ defines
a finite morphism  of
degree ${\rm Vol}(\Delta)$
\[ \Psi \, \, : \P_{\Delta} \to \P^d,\; \;  
p \mapsto (G_0(p): G_1(p): \cdots : G_d(p)). \]
Consider the $\check{\rm C}$ech cocycle
 \[ 
  \alpha = 
  \frac{z_0^d}{z_1 \cdots z_d} \, 
  d \left(\frac{z_1}{z_0} \right) 
  \wedge \cdots \wedge 
  d \left(\frac{z_d}{z_0} \right)
\]
in the standard open covering $\mathcal{U} = \{ z_i \neq 0\}_{i = 0,\ldots,d}$ 
of $\P^d$ with the homogeneous coordinates $z_0,z_1, \ldots, z_d$. One can
show that $\alpha$ determines a generator of $H^d(\P^d,\Omega^d_{\P^d})$ 
\cite[Chapter III, \S~7]{Hartshorne1}. It is easy to check that the map 
\[
 \Psi^*\, : \,H^d(\P^d, \Omega^d_{\P^d}) \to 
              H^d(\P_{\Delta},\Omega^d_{\P_{\Delta}})
\]
sends $\alpha$ to the $\check{\rm C}$ech cocycle
\[
  \Psi^*(\alpha) = 
  \frac{J_G}{G_0\cdots G_d}\,
  \frac{d t_1}{t_1}\wedge \cdots \wedge \frac{d t_d}{t_d}
\]
in the covering $\mathcal{U}' = \{G_i \neq 0\}_{i = 0,\ldots,d}$ of
$\P_\Delta$. Since the sheaf $\Omega_{\P_\Delta}^d$ is dualizing, there exists
a canonical {\it trace map} 
$T_{\P_\Delta}: H^d(\P_\Delta,\Omega_{\P_\Delta}^d)\rightarrow \C$ (see 
\cite[Chapter III, \S~7]{Hartshorne1}). By the property of the trace map 
\[
                   T_{\P_\Delta}([\Psi^*(\alpha)]) = 
  \deg(\Psi)       T_{\P^d}([\alpha]) = 
  {\rm Vol}(\Delta)T_{\P^d}([\alpha])
\]
(see \cite[Chapter III]{Hartshorne2}).
This naturally explains the normalization ${\rm Res}_G(J_G)= 
{\rm Vol}(\Delta)$  in Definition \ref{Tor_res}.

\bigskip

Let us review some  properties of the toric residue. 
\medskip

Choose a regular sequence $G = (G_0,G_1, \ldots, G_d)$ in
$S^1_{\Delta}$ defined by Laurent polynomials $g_0(t),g_1(t),\ldots,g_d(t)$
supported in $\Delta$ as above. It follows from the regularity of $G$ that the
set $V_g$ of common zeros of $G_1,\ldots,G_d$ in  $\P_\Delta$ is 
finite.  Next statement follows immediately from \cite[Theorem~0.4]{CCD} and
\cite[Proposition~1.3]{CDS} and claims that toric residue can be expressed
as a sum of local Grothendieck residues. 

\begin{theo} 
\label{sum_Groth}
Let $p(t_1,\ldots,t_d)$ be a Laurent polynomial with support in the interior
of $(d + 1)\Delta$. We set 
$P := t_0^{d + 1} p(t_1,\ldots,t_d)\in I^{d+1}_{\Delta}$ and choose a 
sufficiently generic
regular sequence 
$(G_0,G_1,\ldots, G_d) = 
(t_0 g_0,t_0 g_1,\ldots, t_0 g_d)\subset S_\Delta^1$ such  the
set $V_g$ of common zeros of $G_1,\ldots,G_d$ in  $\P_\Delta$ is 
contained in $\T \cong (\C^*)^d$. Then 
\[ 
  {\rm Res}_G(P) =  
  \sum_{\xi \in V_g} {\rm res}_{g,\xi} 
  \left( \omega_P \right),
\]
where 
\[
  \omega_P =   
  \frac{p/g_0}{g_1 \cdots g_n} \,
  \frac{d t_1}{t_1} \wedge \cdots \wedge 
  \frac{d t_d}{t_d}
\]
and ${\rm res}_{g,\xi}(\omega_P)$ is the local Grothendieck residue of the
form $\omega_P$ at the point $\xi \in V_g$. 

In particular, if all the common roots of $g_1,\ldots,g_d$ in $\T$ are simple,
then 
\[ 
  {\rm Res}_G(P) =  
  \sum_{\xi \in V_g}
  \frac{ p(\xi)}{g_0(\xi) {J}_g^0(\xi)}, 
\]
where 
\[ 
  J_g^0 :=\det \left( t_i \frac{\partial g_j}
  {\partial t_i} \right)_{1 \leq i, j \leq d}. 
\]
\end{theo} 

\medskip

Let $A$ be a finite subset in $\Delta \cap M$ which contains all vertices
of $\Delta$. Write each of the polynomials $g_0(t),g_1(t),\ldots,g_d(t)$ as
\[
  g_j(t) = \sum_{m\in A} a^{(j)}_m t^m,\quad j = 0,1,\ldots,d.
\]
Let $\Q[a]$ be the polynomial ring in the variables $a^{(j)}_m$ $(m \in A, 
0 \leq j \leq d)$. Denote by
${\mathcal R}_{g_0,\ldots,g_d}(a) \in \Q[a]$ {\it the (unmixed) sparse
$A$-resultant of} $g_0,\ldots,g_d$ defined in \cite[\S 8.2]{GKZ}. The
following statement, which is a reformulation of \cite[Theorem~1.4]{CDS} and
\cite[Proposition~3.5]{Dickenstein}, claims that toric residue is a rational
function in coefficients of Laurent polynomials $g_0,g_1,\ldots,g_d$. 

\begin{theo} 
\label{res_is_rat}
For any interior lattice point $m$ in  $(d+1)\Delta$, 
there exists a polynomial
$Q_m(a)\in \Q[a]$ such that 
\[ 
  {\rm Res}_G(t_0^{d+1} t^m) = 
  \frac{Q_m(a)}{{\mathcal R}_{g_0,\ldots,g_d}(a)}. 
\] 
\end{theo} 

We will be mostly interested in toric residues in the special case when the
regular sequence $F = (F_0,F_1,\ldots,F_d)$ of elements in $S_\Delta^1$ is
constructed as follows. Take a generic Laurent polynomial 
\[ 
  f(t) =  \sum_{m\in A} a_m t^m \in 
  \C[t_1^{\pm 1}, \ldots, t_d^{\pm 1}] 
\]
and define 
\[
  F_0:= t_0 f(t),\, F_1 := t_0 f_1(t),\ldots, F_d := t_0 f_d(t), 
\]
where $f_i(t) := t_i \partial f/ \partial t_i$ $( 1 \leq i \leq d)$. 
In this case, the Jacobians $J_F$, $J^0_F$ become Hessians 
\[ 
  H_f:= \det \left( \left(t_i \frac{\partial}{\partial t_i}\right)
  \left(t_j \frac{\partial}{\partial t_j}\right) t_0 f 
  \right)_{0 \leq i, j \leq d}, \;\;  
  H^0_f:= \det \left(t_j \frac{\partial f_i}{\partial t_j}
               \right)_{1 \leq i, j \leq d}.
\]
We write for simplicity $S_f$, $I_f$, ${\rm Res}_f$ instead of $S_F$, $I_F$, 
${\rm Res}_F$ respectively. 

The {\it principal $A$-determinant} $E_A(f)$ of polynomial $f(t)$
is defined to be the sparse resultant ${{\mathcal R}_{f,f_1,\ldots,f_d}(a)}$
of polynomials $f(t), f_1(t), \ldots, f_d(t)$ \cite[\S 10.1]{GKZ}. 
It follows from \cite[\S 4]{Batyrev1} that principal $A$-determinant
$E_A(f)$ is nonzero if and only if $F_0,F_1,\ldots,F_d$ form 
a regular sequence in
$S_\Delta^1$. In the latter case the polynomial $f(t)$ will be called 
{\it $\Delta$-regular}. Note that  $\Delta$-regularity of $f$ guarantees that
all critical points of $f$ are isolated. 

\begin{rem} There exists another definition of $\Delta$-regularity of a
Laurent polynomial $f = \sum_{m \in A} a_m t^m$. For any face 
$\Gamma$ of $\Delta$, we define the Laurent polynomial 
\[
  f_\Gamma := \sum_{m \in A \cap \Gamma} a_m t^m.
\]
The polynomial $f$ is called $\Delta$-{\it regular} if 
\begin{equation}
\left \{ t \in \T \; : \; f_\Gamma(t) = \frac{\partial f_\Gamma}{\partial t_1}
 (t) = \cdots = \frac{\partial f_\Gamma}{\partial t_d}
 (t) = 0 \right\} = 
\emptyset 
\label{DR}
\end{equation} 
for all faces $\Gamma \subseteq \Delta$. This definition is equivalent to the
previous one, because  $E_A(f)$ is equal to the product of
discriminants \cite[Chapter 10, Theorem 1.2]{GKZ} 
\[ 
  E_A(f) = \pm \prod_{\Gamma \subseteq \Delta} D_{A\cap\Gamma} 
  (f_{\Gamma})^{\mu_\Gamma}, 
  \;\;\mu_\Gamma \in \Z_{> 0}
\]
and for each  face $\Gamma \subseteq \Delta$  the  condition (\ref{DR}) holds  
if and only if the discriminant $D_{A\cap\Gamma}(f_{\Gamma})$ is nonzero. 
\label{discrim}
\end{rem} 

\medskip

The statements of theorems \ref{sum_Groth} and \ref{res_is_rat} for 
the toric residue ${\rm Res}_f$ can be summarized as follows: 

\begin{theo} 
\label{Prop_tor_res}
Let $P = t_0^{d + 1} p(t_1,\ldots,t_d)$ be an arbitrary element in
$I_{\Delta}^{d + 1}$. Then 

$(i)$ for any sufficiently generic   $\Delta$-regular Laurent 
polynomial $f(t)$ such that the set 
$V_f$ of all common zeros of $F_1, \ldots, F_d$ in $\P_\Delta$ is contained in 
$\T \cong (\C^*)^d$ and any critical point $\xi \in V_f$ of $f$ is nondegenerate
(i.e., $H_f^0(\xi) \neq 0$), one has 
\[ 
  {\rm Res}_{f}(P) =  \sum_{\xi \in V_f}
  \frac{p(\xi)}{f(\xi) {H}_f^0(\xi)}.  
\]

$(ii)$ the toric residue ${\rm Res}_f(P)$ is a rational function 
in the coefficients of the polynomials $p(t)$ and $f(t)$. 
In particular, 
for any interior lattice point $m$ in  $(d+1)\Delta$, 
there exists a polynomial
$Q_m(a)\in \Q[a]$ such that 
\[ 
  {\rm Res}_f(t_0^{d+1} t^m) = 
  \frac{Q_m(a)}{E_A(f)}. 
\] 
\end{theo}  

\medskip

Let $N:= Hom(M,\Z)$ be the dual to $M$ lattice and 
$N_{\R} := N \otimes \R$. We denote by $\< *, * \>$ the natural pairing
$M_{\R} \times N_{\R} \to \R$. 

\begin{dfn}[\cite{Batyrev2}]
A polytope $\Delta\subset M_\R$ with vertices in $M$ is called 
{\it reflexive} if it contains $0$ as interior point and its polar polytope 
\[
  \Delta^* = \{ y \in N_\R:\<x,y \>\ge -1, \quad \forall x \in 
  \Delta\}\subset N_\R
\]
has vertices in $N$. We remark that  $\Delta$ is reflexive if and only if 
$\P_\Delta$ is a Gorenstein toric Fano variety and ${\mathcal O}_{\P_\Delta}(1)$ 
is the anticanonical sheaf on $\P_\Delta$.
\end{dfn}

The reflexivity of $\Delta$ implies that $\Delta^*$ is also reflexive and
$(\Delta^*)^* = \Delta$. We will demand that the finite subset 
$A \subset M \cap \Delta$ contains not only all the vertices of $\Delta$, but
also the lattice point $0$ as interior. 

\medskip

Now consider the toric residue ${\rm Res}_f$ in the special case when the
support polytope $\Delta$ of $f$ is reflexive. In this case,  $I_{\Delta}$ is
the principal ideal 
in $S_{\Delta}$ generated by $t_0$. This implies that 
$H_f \in I^{d+1}_{\Delta} \subset S_{\Delta}^{d+1}$ can be uniquely written as
product $t_0 H'_f$, where  
\[
  H'_f = t^d_0 \det
  \begin{pmatrix} 
  f &  f_1 &  \cdots & f_d  \\
  f_1 & t_1 \frac{\partial f_1}{\partial t_1} & \cdots &   
  t_1 \frac{\partial f_d}{\partial t_1} \\
  \vdots &  \vdots  &  \ddots  & \vdots  \\
  f_d & t_d \frac{\partial f_1}{\partial t_d} & \cdots &   
  t_d \frac{\partial f_d}{\partial t_d} 
  \end{pmatrix}
  \in S^{d}_{\Delta}.
\]
The ring $S_{\Delta}$ is Gorenstein and the multiplication in Artinian graded
ring $S_f$ defines the perfect pairings 
\[
  S_f^k \times S^{d - k}_f \to S^d_f, \quad k = 0,\ldots,d,
\]
where $S^d_f$ is a $1$-dimensional $\C$-vector space generated by the class of
$H'_f$.
By abuse of notations, we denote by ${\rm Res}_f$ also the $\C$-linear map 
\[
  {\rm Res}_f : S_\Delta^d \rightarrow \C, 
\]
vanishing on the subspace  $\< F_0, F_1, \ldots, F_d\>S_\Delta^{d-1} 
\subset  S_\Delta^d$ such that  ${\rm Res}_f(H'_f)= 
{\rm Vol}(\Delta)$. This map induces a canonical isomorphism 
\begin{equation}
\label{Res_f_T}
 S^d_f =  S_\Delta^d /\<F_0,F_1,\ldots,F_d\> S_\Delta^{d-1} \cong \C.
\end{equation}

\bigskip

\section{Morrison-Plesser moduli spaces}

Let $\P = \P_\Si$ be a $d$-dimensional projective simplicial toric variety
over $\C$ defined by a simplicial fan $\Sigma$ whose $1$-dimensional cones are
generated by integral vectors $e_1, \ldots, e_n \in N$. There exists a
canonical embedding $M = Hom(N, \Z)\hookrightarrow \Z^n$ defined by 
\[ 
  m \mapsto (\< m, e_1\>, \ldots, \<m, e_n\>) 
\]
which gives rise to the short exact sequence 
\[ 
  0 \rightarrow M \rightarrow \Z^n \rightarrow 
  {\rm Cl}(\P) \rightarrow 0, 
\]
where  ${\rm Cl}(\P)$ is canonically isomorphic to the group of Weil 
 divisor classes on $\P$ modulo linear equivalence.

Let $G \subset (\C^*)^n$ be the diagonalizable algebraic group over $\C$ with
the character group ${\rm Cl}(\P)$. Then $\P$ can be described   as a quotient
$U(\Sigma)/G$ (space of orbits), where $U(\Sigma)$ is an  open dense subset
in $\C^n$ and the action of $G$ on $U(\Sigma)$ is induced by the embedding 
$G \hookrightarrow (\C^*)^n$ defined by the epimorphism $\Z^n \rightarrow 
  {\rm Cl}(\P)$. The standard affine coordinates 
$z_1, \ldots,z_n$ on $\C^n$ determine {\em homogeneous coordinates} on $\P$
\cite{Cox1}. Every equation $z_j = 0$ $( 1 \leq j \leq n)$ defines a Weil
($\Q$-Cartier) divisor $D_j$ on $\P$. We denote by  $\chi_1, \ldots, \chi_n$
the characters of $G$ defining the representation of $G$ in $GL(n,\C)$. These
characters correspond to the Weil divisor classes 
$[D_1], \ldots, [D_n] \in {\rm Cl}(\P)$.
Let  $R(\Si)$ be the subgroup in $\Z^n$ consisting of all lattice vectors
$\lambda = (l_1,\ldots,l_n)$ such that $l_1 e_1 + \cdots + l_n e_n = 0$ and 
 $N'$ the sublattice of finite index in $N$ generated by 
$e_1, \ldots, e_n$. We have the exact sequence 
\[  
  0 \rightarrow R(\Sigma) \rightarrow \Z^n \rightarrow 
  N  \rightarrow N/N' \to 0, 
\]
where the middle map is $(l_1, \ldots, l_n) \mapsto \sum_{j =1}^n l_j e_j$. 
Tensoring  this sequence by $\Q$, one gets 
\[  
  0 \rightarrow R(\Sigma) \otimes \Q  \rightarrow \Q^n \rightarrow 
  N \otimes \Q  \rightarrow 0. 
\]
Comparing the latter with the exact sequence 
\[ 
  0 \rightarrow M \otimes \Q \rightarrow \Q^n \rightarrow 
  {\rm Cl}(\P) \otimes \Q \rightarrow 0, 
\]
we obtain the canonical isomorphisms 
\[ 
  Hom(R(\Sigma), \Q) \cong  {\rm Cl}(\P) \otimes \Q \cong  
  {\rm Pic}(\P) \otimes \Q \cong H^2(\P, \Q). 
\]
Therefore, $R(\Sigma)_\Q:=R(\Sigma) \otimes \Q$ can be identified with the
rational homology group $ H_2(\P, \Q)$. For any $\Q$-divisor 
$D = \sum_{j =1}^n c_j D_j $ $( c_j \in \Q)$ and for any 
$\lambda  = (l_1, \ldots, l_n) \in R(\Sigma)_\Q$, one has the intersection
number 
\[ 
  (D, \lambda) = \sum_{i =1}^{n} c_i l_i \in \Q. 
\]
If $H$ is an ample Cartier divisor on $\P$, then the fan $\Sigma$ can be
obtained as a normal fan for the simple $d$-dimensional polytope 
\[
  \Delta_H := \left\{ (x_1, \ldots, x_n) \in \R^n_{\geq 0} \; : \; 
  \sum_{j=1}^n (D_j, \lambda)x_j = (H, \lambda), \quad \forall \lambda \in  
  R(\Sigma)_{\Q}\right\}.
\]
All vertices of $\Delta_H$ belong to the sublattice $M_H \subset \R^n$, where
$M_H$ is the set of all integral vectors $(x_1, \ldots, x_n) \in \Z^n$
such that $H$ is linearly equivalent to $\sum_{j =1}^n x_j D_j$. 

Let $K_{\rm amp}(\P)$ be the closed {\it ample $($or K\"ahler$)$ cone} in
${\rm Cl}(\P) \otimes \R \cong H^2(\P, \R)$ and $K_{\rm eff}(\P)$ the dual
to $K_{\rm amp}(\P)$ {\it Mori cone} of effective curves in 
$R(\Si)_{\R} \cong H_2(\P, \R)$. The cone of vectors 
$\beta = (b_1, \ldots, b_n)  \in  R(\Si)_{\R} = R(\Sigma) \otimes \R$ 
such that $b_1, \ldots, b_n \ge 0$ will be denoted by $K^+(\P)$. It is easy
to see that $K^+(\P)$ is always a subcone of the Mori cone 
$K_{\rm eff}(\P)$. 

For any lattice point $\beta=(b_1, \ldots, b_n) \in K_{\rm eff}(\P)$, we
will construct a simplicial toric variety $\P_{\beta}$ which can be
considered as a ``naive compactification'' of the moduli space of rational
maps $\phi\;  : \; \P^1 \to \P$  such that the class 
$[\phi(\P^1)] \in H_2(\P,\R)$ is equal to $\beta$. In the case $\beta = 0$, 
the toric variety  $\P_{\beta}$ coincides with $\P$.  

First we consider the case  $\beta = (b_1, \ldots, b_n) \in K^+(\P)$,
i.e., all the $b_1, \ldots, b_n$ are supposed to be nonnegative. Since 
$(D_j, \beta) = b_j$ $(1\leq j \leq n)$, we can construct a map 
$\phi\, : \, \P^1 \to \P$ with $[\phi(\P^1)] = \beta$ by choosing $n$
homogeneous binary forms $\phi_1(u,v), \ldots, \phi_n(u,v) \in \C[u,v]$ such
that ${\rm deg}\, \phi_j = b_j$ $( 1\leq j \leq n)$. Let  $ \C_j(\beta)$ be
the  space of all homogeneous binary forms of degree $b_j$ 
$({\rm dim}\,  \C_j(\beta) = b_j + 1$). We set 
$\C(\beta): = \bigoplus_{j =1}^n \C_j(\beta)$  and denote by 
$z_0^{(j)},\ldots, z_{b_j}^{(j)}$ the coordinates on $\C_j(\beta)$
corresponding to the standard monomial basis of  $\C_j(\beta)$. We define the
action of $G$ on   $\C_j(\beta)$ as the scalar multiplication 
by the character
$\chi_j$ $( 1\leq j \leq n)$. This defines an effective action of $G$ on
$\C(\beta)$. If $\chi_H = \chi_1^{c_1}\cdots \chi_n^{c_n}$ is the character
of $G$ corresponding to the class of the ample Cartier divisor 
$H = \sum_{j = 1}^n c_j D_j$ on $\P$, then we define the 
{\bf Morrison-Plesser moduli space} $\P_\beta$ as the GIT-quotient
$\C(\beta)//G$ with respect to 
the linearization by $\chi_H$ of the structure sheaf on  $\C(\beta)$. We
define the dense open subset $U(\beta) \subset \C(\beta)$ as the union 
of all open
subsets 
\[ 
  U_{i_1, \ldots, i_{n-d}}(\sigma) = 
  \left\{ z \in \C(\beta) \; : \; z_{i_1}^{(j_1)} z_{i_2}^{(j_2)} \cdots 
  z_{i_{n-d}}^{(j_{n-d})} \neq 0 \right\}, 
\]
where $\sigma$ runs over all $d$-dimensional cones of $\Sigma$, 
$\{ e_{j_1}, e_{j_2}, \ldots, e_{j_{n-d}}\}$ is the set of all vectors from
$\{e_1, \ldots, e_n \}$ which do not belong to $\sigma$, 
and each index $i_k$ 
$(1 \leq k \leq n - d)$ runs independently over all elements of 
$\{0,1,\ldots, b_k\}$. It is easy to show that the Morrison-Plesser moduli
space $\P_\beta$ is also 
the space of orbits $U({\beta})/G$. Moreover, $\P_\beta$
is a projective simplicial toric variety of dimension 
$d + \sum_{j =1}^n b_j$. 

In general, the cone $K^+(\P)$ is smaller than $K_{\rm eff}(\P)$. In
Section~\ref{Sect_Hirz} we consider  such a situation for  
$\P = \mathbb{F}_1$ and  show that there exist
infinitely many   classes $\beta \in K_{\rm eff}(\mathbb{F}_1) \setminus 
K^+(\mathbb{F}_1)$ which can not be represented by irreducible curves 
$C \subset  \mathbb{F}_1$ (see Remark~\ref{irred-cl}). Therefore, if one of
the coordinates $b_j$ of $\beta = (b_1, \ldots, b_n)$ is negative, it may happen
that there is no rational map $\phi\, : \, \P^1 \to \P$ such that 
$[ \phi(\P^1)] = \beta$, but the  corresponding  
Morrison-Plesser moduli space
$\P_\beta$ (see Definition \ref{DMP} below) is not empty.

Now let $\beta = (b_1, \ldots, b_n)$ be an arbitrary lattice point in 
$K_{\rm eff}(\P)$. For any $j \in \{ 1, \ldots, n\}$, we define the free
abelian group $\Z_j(\beta)$ as 
\[ 
  \Z_j(\beta) := \left\{ \begin{array}{ll} \Z^{b_j+1}, & 
  \mbox{\rm if $b_j \geq 0$},\\
  0, &  \mbox{\rm if $b_j < 0$}. \end{array} \right.    
\]
Using the standard basis of $\Z_j(\beta)$, we write each element of
$\Z_j(\beta)$ as the integral vector 
$(x_0^{(j)}, x_1^{(j)}, \ldots, x_{b_j}^{(j)})$. We set 
$\Z(\beta) := \bigoplus_{j =1}^n \Z_j(\beta)$, 
$\R(\beta) := \Z(\beta) \otimes \R$ and 
denote by $\R_{\geq 0}(\beta)$ the set
of all vectors in $\R(\beta)$ having nonnegative coordinates. 

\begin{dfn} 
Let be $r$ the number of negative coordinates of a lattice point 
$\beta= (b_1,\ldots,b_n) \in K_{\rm eff}(\P)$. Without loss of generality
we may assume that $b_1,\ldots, b_{n-r} \geq 0$ and $b_{n-r+1}, \ldots, 
b_n < 0$. We
define the convex set 
\[ 
  \Delta_H^\beta := \left\{ x \in \R_{\geq 0 }(\beta) \; : \; 
  \sum_{j=1 }^{n-r} (D_j, \lambda) \left( \sum_{i =0}^{b_j} 
  x_i^{(j)} \right)  = (H, \lambda), \quad \forall \lambda \in  
  R(\Sigma)_{\Q}\right\}. 
\]
Denote by  $M_H^{\beta} \subset \Z(\beta)$ the set of all lattice points 
$x \in \Z(\beta)$ whose coordinates $x_i^{(j)}$ satisfy the condition: the
divisor 
\[ 
  \sum_{j = 1 }^{n-r} \left( \sum_{i =0}^{b_j} 
  x_i^{(j)} \right)D_j  
\] 
is linearly equivalent to $H$ up to a linear combination of $D_{n-r+1}, 
\ldots, D_n$. 
\end{dfn} 

\begin{prop} 
$\Delta_H^\beta$ is a compact convex simple  polytope having  
vertices in
$M_H^{\beta}$. The normal fan $\Sigma_\beta$ of the polytope $\Delta_H^\beta$
does not depend on the choice of the ample divisor $H$. 
\end{prop} 

\begin{proof} In order to see the compactness of $\Delta_H^\beta$, we remark
that there exists a canonical affine linear mapping 
$\pi_\beta : \Delta_H^\beta \to \Delta_H^0 = \Delta_H$ which replaces every
$b_j + 1$ coordinates $x_0^{(j)}, x_1^{(j)}, \ldots, x_{b_j}^{(j)}$ by
their sum $x_j = \sum_{i =0}^{b_j} x_i^{(j)}$. It is clear that 
$\pi_\beta(\Delta_H^\beta) \subset 
\R^n$ is a face of $\Delta_H^0$ defined by the
equations $x_j = 0$ $(n-r \leq j \leq n)$. Therefore  
${\rm dim}\, \pi_\beta(\Delta_H^\beta) = d-r$ if 
$\pi_\beta(\Delta_H^\beta)$ is not empty (the latter holds if and only
if $e_{n-r}, \ldots, e_n$ generate a $r$-dimensional cone in $\Sigma$).
We observe that the preimage
$\pi_\beta^{-1}(x)$ of a point in $x \in \pi_\beta(\Delta_H^\beta)$ is the
product of $n-r$ simplices of dimensions $b_1, \ldots, b_{n-r}$. 
Therefore, $\Delta_H^\beta$ is compact 
and  
$${\rm dim}\, \Delta_H^\beta = d-r +   
\sum_{j=1}^{n-r} b_j.$$
Now we want to describe all faces of codimension $1$ (i.e., facets) of 
$\Delta_H^\beta$. 
It is clear that each facet must be defined  by an equation 
$x_i^{(j)} = 0$ for some $0 \leq i \leq b_j$, $1 \leq j \leq n-r$. 
However, it is not true 
in general that every equation 
$x_i^{(j)} = 0$ defines a facet. It is easy to show that 
the equation 
$x_i^{(j)} = 0$ defines a facet of $\Delta_H^\beta$ if and only if 
either $b_j >0$, or if $b_j =0$ and $e_j, e_{n-r}, \ldots, e_n$ generate a 
$(r+1)$-dimensional cone in $\Sigma$. The $\pi_\beta$-image of a vertex 
$y \in \Delta_H^\beta$ is a vertex of $\Delta_H$. Take an arbitrary 
vertex $x \in \pi_\beta(\Delta_H^\beta)$. 
Since  $x$ is a vertex of $\Delta_H$,  
there exists a subset $\{{i_1}, \ldots, {i_d} \} \subset \{1, \ldots, n\}$ 
which contains $\{ n-r+1, \ldots, n\}$ such that 
$e_{i_1}, \ldots, e_{i_d}$ are generators of a $d$-dimensional cone 
in $\Sigma$ and $x \in \Delta_H$ is defined by the conditions
$x_{i_1} =  \cdots =  x_{i_d} = 0$. Since $\Delta_H$ is a simple polytope
we have $x_j > 0$ for all $j \not\in \{{i_1}, \ldots, {i_d} \}$ and the 
equation 
$x_j = \sum_{i =0}^{b_j} x_i^{(j)}$ defines a $b_j$-dimensional simplex 
with  $(b_j+1)$-vertices. In this way, 
we obtain $\prod_j (b_j +1)$ ($j$ runs over 
$\{ 1, \ldots, n\} \setminus \{{i_1}, \ldots, {i_d} \}$) vertices 
of $\Delta_H^\beta$ as $\pi_\beta$-preimages of $x$. By this method, 
we get  
all vertices of  $\Delta_H^\beta$ as $\pi_\beta$-preimages of vertices 
of $\pi_\beta(\Delta_H^\beta)$. Moreover, one sees that each vertex of 
$\Delta_H^\beta$ is contained in exactly ${\rm dim}\, \Delta_H^\beta$
facets: for each $j \in \{{i_1}, \ldots, {i_d} \} \setminus 
\{ n-r+1, \ldots, n\}$ we get $b_j+1$ facets $x_j^{(i)} =0$ 
$(0 \leq i \leq b_j)$ containing a chosen $\pi_\beta$-preimage $y$ of 
$x$ and for  each $j \in \{1, \ldots, n\} \setminus  
\{{i_1}, \ldots, {i_d} \}
$ we get $b_j$ facets $x_j^{(i)} =0$ containing $y$. 
Therefore $\Delta_H^\beta$ is a simple polytope. Since the combinatorial
structure of  $\Delta_H^\beta$ is completely determined by $\beta$ and 
$\Delta_H$,  it  does not depend on the  choice of $H$ and the same 
is true for the normal fan $\Sigma_\beta$.
\end{proof}

\begin{dfn} 
Let $\beta = (b_1, \ldots, b_n)$ be an arbitrary  lattice point in 
$K_{\rm eff}(\P)$. The projective simplicial toric variety $\P_\beta$
associated with the normal fan $\Sigma_\beta$ of the polytope $\Delta_H^\beta$
is called the {\bf Morrison-Plesser moduli space} corresponding to 
$\beta \in K_{\rm eff}(\P)$. 
\end{dfn} 

It follows from the proof of the last proposition that 
\[
  {\rm dim} \, \P_\beta = {\rm dim}\, \Delta_H^\beta =
 d -n + \sum_{j=1}^{n-r} (b_j + 1)
\]
if $\P_\beta$ is not empty.   
It is easy to see that the last definition of the Morrison-Plesser moduli
space coincides with the previous one in the case $\beta \in K^+(\P)$. We also 
remark that $\P_\beta$ is 
nonsingular for all $\beta \in K_{\rm eff}(\P)$ if  $\P$ is nonsingular. 
We will need the following property of the Morrison-Plesser 
moduli spaces $\P_\beta$: 

\begin{prop} 
\label{DMP}
There exists a canonical surjective homomorphism 
$$\psi_\beta : H^2(\P, \Q) \to H^2(\P_{\beta}, \Q)$$ which is always 
bijective  if $\beta \in K^+(\P)$. 
\end{prop} 

\begin{proof} 
Let $\{ e_{j_1}, \ldots, e_{j_k}\}$ be the set of all generators 
$e_j \in \{ e_1, \ldots, e_{n-r} \}$ 
such that $b_j =0$ and $e_j, e_{n-r+1}, \ldots, e_n$ do not generate 
a $(r+1)$-dimensional cone in $\Sigma$. We 
define the subgroup $G_\beta \subset G$ to be  the common 
kernel of the  characters $\chi_{j_1}, \ldots, \chi_{j_k}$, i.e.,
\[
  G_\beta := \{ g \in G \: : \; 
  \chi_{j_1}(g) = \cdots = \chi_{j_k}(g) =1 \}.
\] 
Then the simplicial toric 
variety $\P_\beta$ can be obtained as a geometric quotient of 
an affine space of dimension $(b_1 +1) + \cdots + (b_{n-r}+1) -k$ modulo 
the linear action of $G_\beta$. The embedding 
$G_\beta \hookrightarrow G$ induces the surjective homomorphism 
of the character groups ${\rm Cl}(\P) \to {\rm Cl}(\P_\beta)$. This 
homomorphism is bijective  if $\beta \in K^+(\P)$, because 
in the latter case $G = G_\beta$. Tensoring by $\Q$, we obtain 
the canonical surjective homomorphism 
$\psi_\beta : H^2(\P, \Q) \to H^2(\P_{\beta}, \Q)$. 
\end{proof}

\begin{dfn}
Assume that the anticanonical class $-K_{\P}$ 
of $\P$ is nef, i.e., $(-K_{\P}, \beta) = 
\sum_{j =1}^n b_j \geq 0$ for 
all $\beta \in K_{\rm eff}(\P)$ and $r$ is the number of negative coordinates
of $\beta$, i.e., 
$b_1,\ldots, b_{n-r} \geq 0$ and $b_{n-r+1}, \ldots, b_n < 0$.
 By abuse of notations, let us denote by $[D_j] \in H^2(\P_\beta, \Q)$ 
($1 \leq j \leq n$) also the image of $[D_j] \in H^2(\P, \Q)$ under
$\psi_\beta$. Using the multiplication in the cohomology ring 
$H^*(\P_\beta, \Q)$, we define  the intersection product 
\[ 
  \Phi_{\beta} := 
  ([D_1] + \cdots + [D_n])^{b_1 + \cdots + b_n} 
  \prod_{j=n-r+1}^n [D_j]^{-b_j - 1},
\]
considered as a cohomology class in $H^{2(\dim \P_\beta -d)}(\P_{\beta}, \Q)$
and call $\Phi_{\beta}$ the {\bf Morrison-Plesser class} of $\P_\beta$. 
\end{dfn}

\bigskip

\section{Toric Residue Mirror Conjecture}

In order to formulate our conjecture, we need some results about Newton
polytopes of principal $A$-determinants due to Gelfand, Kapranov and
Zelevinsky \cite{GKZ}. 
\medskip

Let $\Delta \subset M_\R$ be a $d$-dimensional polytope with vertices in
$M$. Denote by $A$ a finite subset in $\Delta \cap M$ which includes all
vertices of $\Delta$. 

\begin{dfn}
By a {\it triangulation  ${\mathcal T} = \{\tau_1, \ldots, \tau_k\}$ of
$\Delta$ associated with $A$}, we mean a decomposition of $\Delta$ into a
union of $d$-dimensional  simplices $\tau_1, \ldots, \tau_k$ having 
 vertices in
$A$ such that any nonempty  intersection $\tau_i \cap \tau_j$ 
is a common face
of $\tau_i$ and $\tau_j$. A triangulation ${\mathcal T}$ associated with $A$
is called {\it coherent} if there exists a convex piecewise-linear function 
$\phi : \Delta \to \R$ whose domains of linearity are precisely the simplices of
${\mathcal T}$. 
\end{dfn} 

\begin{dfn}
Denote by $\R^A$ the space of all real-valued functions on $A$. Let 
${\mathcal T}$ be  a triangulation of $\Delta$ associated with $A$. The
function $\chi_{\mathcal T}: A\rightarrow \R$ defined as 
\[
  m \mapsto \sum_{i\, :\,  m\in{\rm Vert}(\tau_i)}{\rm Vol}(\tau_i),
\]
where the sum of the normalized volumes ${\rm Vol}(\tau_i)$ runs over all
simplices of $\tau_i \in {\mathcal T}$ containing $m \in A$ as vertex, is
called the {\it characteristic function} of ${\mathcal T}$. The 
{\bf secondary polytope} ${\rm Sec}(A)$ is defined as the convex hull of
the vectors $\chi_{\mathcal T} \in \R^A$, where  ${\mathcal T}$ runs over all
triangulations of $\Delta$ associated with $A$. 
\end{dfn}

\begin{theo}[\cite{GKZ}, Chapter~7.1] 
The secondary polytope ${\rm Sec}(A) \subset \R^A$ is a 
$(|A| - d - 1)$-dimensional polytope whose vertices are exactly the
characteristic functions $\chi_{\mathcal T}$ corresponding to all coherent
triangulations ${\mathcal T}$ of $\Delta$. 
\end{theo}

Consider a generic Laurent polynomial 
\[ 
  f(t) = \sum_{m \in A}a_m t^m \in \C[t_1^{\pm 1},\ldots,t_d^{\pm 1}].
\]
The principal $A$-determinant $E_{A}(f)$ is a certain polynomial in 
$|A|$ independent variables $\{ a_m\}_{m \in A}$  with
integral coefficients. The following theorem will be  very important in the
sequel.

\begin{theo}[\cite{GKZ}, Chapter~10.1, Theorem 1.4] 
\label{E_A=Sec(A)} 
The Newton polytope of $E_A(f)$ coincides with the secondary polytope 
${\rm Sec}(A)$. If ${\mathcal T} = \{\tau_1,\ldots,\tau_k\}$ is a coherent
triangulation corresponding to some vertex of ${\rm Sec}(A)$, then the
coefficient at the monomial $\prod_{m \in A} a_m^{\chi_{\mathcal T}(m)}$ in
$E_A(f)$ is equal (up to sign) to  the product 
\[
   \prod_{i =1}^k {\rm Vol}(\tau_i)^{{\rm Vol}(\tau_i)}.
\]
\end{theo}

Recall the notion of a Laurent series of a rational function at a vertex of the
Newton polytope of its denominator (see, e.g., \cite{GKh} or
\cite[p.~195]{GKZ}). 

\begin{dfn} 
Let $P(a), Q(a) \in \C[a_1^{\pm 1},\ldots,a_n^{\pm 1}]$ be two arbitrary
Laurent polynomials and let $v \in \Z^n$ be a vertex of the Newton polytope of 
$Q = \sum_{w} c_w a^w$. We write $Q(a) = c_v a^v (1 + \tilde{Q}(a))$, where 
\[ 
  \tilde{Q}(a) := \sum_{w \ne v} \frac{c_{w}}{c_v} a^{w - v}. 
\]
It is easy to see that each Laurent monomial in $a_1, \ldots, a_n$ appears
with nonzero coefficient in $(\tilde{Q}(a))^i$ only for finitely many values of
$i$. So the expression 
\[ 
  \frac{1}{1 + \tilde{Q}(a)} = \sum_{i = 0}^{\infty} (-1)^i (\tilde{Q}(a))^i 
\]
is well-defined as a Laurent power series in the variables 
$a_1, \ldots, a_n$. The product 
\[
  P(a)\cdot c_v^{-1} a^{-v} \cdot 
  (1 - \tilde{Q}(a) + (\tilde{Q}(a))^2 - \cdots)
\]
is called the {\bf Laurent series of the rational function  $P(a)/Q(a)$ at the
vertex $v$} of the Newton polytope of $Q$. 
\end{dfn}

Now we are able to  formulate our Toric Residue Mirror  Conjecture: 

\begin{conj} 
\label{TRMC}
Let $\Delta \subset M_{\R}$ be an arbitrary  reflexive $d$-dimensional
polytope and $A$ a finite subset in $\Delta \cap M$ containing $0$ and all
vertices of $\Delta$. Choose any  coherent triangulation 
${\mathcal T} =\{\tau_1, \ldots, \tau_k\} $ of $\Delta$ associated with $A$
such that $0$ is a vertex of all the simplices 
$\tau_1, \ldots, \tau_k$. Denote by $\P = \P_{\Si({\mathcal T})}$ the simplicial
toric variety defined by the fan $\Sigma = \Sigma({\mathcal T}) \subset M_\R$
whose $d$-dimensional cones are exactly $\sigma_i := \R_{\geq 0} \tau_i$ 
$( 1 \leq i \leq k)$. If $A = \{v_0 = 0, v_1, \ldots, v_n\}$ and 
\[ 
  f(t) :=   1 - \sum_{i = 1}^n a_i t^{v_i},
\]
then for any homogeneous polynomial $P(x_1,\ldots,x_n)\in\Q[x_1,\ldots,x_n]$
of degree $d$ the Laurent expansion of the toric residue 
\[
  R_P(a) := (-1)^d\,{\rm Res}_f(t_0^d\, P(a_1 t^{v_1},\ldots,a_n t^{v_n}))
\]
at the vertex $v_{\mathcal T} \in {\rm Sec}(A)$ corresponding to the coherent
triangulations ${\mathcal T}$ coincides with the generating function of 
intersection numbers 
\[
  I_P(a) := \sum_{\beta\in K_{\rm eff}(\P)} I(P, \beta) a^\beta, 
\]
where the sum runs over all integral points $\beta = (b_1,\ldots,b_n)$ of the
Mori cone $K_{\rm eff}(\P)$, $a^\beta := a_1^{b_1}\cdots a_n^{b_n}$, 
\[
  I(P, \beta) = \int_{\P_\beta}P([D_1],\ldots,[D_n])\Phi_\beta = 
  \< P([D_1],\ldots,[D_n])\Phi_\beta \>_\beta, 
\]
and $\Phi_\beta \in H^{2(\dim \P_\beta -d)}(\P_{\beta}, \Q)$ is 
the Morrison-Plesser class of $\P_{\beta}$. We assume $I(P,\beta)$ 
to be zero if $\P_\beta$ is empty. 
\end{conj}

\begin{rem}
Consider any coherent triangulation 
${\mathcal T} =\{\tau_1, \ldots, \tau_k\} $ of $\Delta$ associated with $A$
such that $0$ is a vertex of all simplices and let $v_{\mathcal T}$ be the
corresponding vertex of the $(n-d)$-dimensional polytope ${\rm Sec}(A)$ as
above. It is easy to show that the cone 
$\R_{\geq 0}({\rm Sec}(A) - v_{\mathcal T})$ can be canonically 
identified with the Mori cone $K_{\rm eff}(\P) \subset H_2(\P,\R)$.  
\end{rem}

We want to check the statement of our conjecture for the coefficient $I(P,0)$
of the power series $I_P(a)$. For this purpose, we choose a $\Sigma$-piecewise
linear function $\varphi$ on $M_\R$ corresponding to an ample Cartier 
divisor $H$ on $\P_\Si$ and make the substitution $a_i = u^{\varphi(v_i)}$ 
($1 \leq i \leq n$), where $u$ is a variable. This substitution 
is equivalent to the
consideration of the $1$-parameter family of Laurent polynomials 
\[ 
  f(t) = 1 - \sum_{i = 1}^n u^{\varphi(v_i)} t^{v_i} 
\]
depending on $u$. We remark that our substitution transforms $I_P(a)$ into a
formal power series 
\[ 
  I_P(u) = \sum_{\beta\in K_{\rm eff}(\P_\Si)} 
  I(P, \beta) u^{(H,\beta)} \in \Q[[u]]. 
\]
Since $(H,\beta)>0$ for all nonzero 
$\beta \in  K_{\rm eff}(\P_\Si) \cap R(\Sigma)$, we have 
\[  
  \lim_{u \to 0}  I_P(u) = I(P, 0) = \int_{\P} P([D_1], \ldots, [D_n]). 
\]

Conjecture \ref{TRMC} for $\beta = 0$ is equivalent to the 
following: 

\begin{theo} 
For any homogeneous polynomial $P(x_1, \ldots, x_n)$ of degree $d$, one has 
\[ 
 (-1)^d \lim_{u \to 0} {\rm Res}_f(t_0^d P(u^{\varphi(v_1)}t^{v_1}, 
  \ldots,u^{\varphi(v_n)}t^{v_n}) = \int_{\P} P([D_1], \ldots, [D_n]).
\]
\end{theo} 

\begin{proof}
Let $\varphi$ be a $\Sigma$-piecewise linear function on $M_{\R}$
corresponding to an ample Cartier divisor on $\P_\Si$ as above. 
Without loss of generality we may assume that $\varphi$ is positively defined,
i.e., $\varphi\geq 0$ on $M_\R$ and $\varphi(x) = 0$ for some $x\in M_\R$ if
and only if $x = 0$. Let
$S_\Delta[u]$ be $\C[u]\otimes_\C S_{\Delta}$, which is considered as a graded
algebra over the polynomial ring $\C[u]$. 
We denote by 
\[
  S_{\varphi}[u] = \bigoplus_{l = 0}^{\infty} S_{\varphi}^l[u] 
\]
the graded $\C[u]$-subalgebra in $S_\Delta[u]$ whose $l$-th homogeneous
component $S_{\varphi}^l[u]$ is spanned as $\C$-vector space by all monomials
$u^r t_0^l t^m$ such that the lattice point $m$ is contained in $l\Delta$ and
$r\ge \varphi(m)$. 
It is easy to see that the set $S_{\varphi}[u]$ is closed under the
multiplication: if $u^r t_0^l t^{m}, u^{r'} t_0^{l'} t^{m'} \in
S_{\varphi}[u]$, then $u^{r + r'}t_0^{l + l'} t^{m + m'}\in
S_{\varphi}[u]$, because $r + r' \ge \varphi(m) + \varphi(m') \ge \varphi(m
+ m')$. 

Let us set $y_0 := t_0$ and  $y_i := - u^{\varphi(v_i)} t_0 t^{v_i}$
($1 \le i \le n$). By definition, the elements $u, y_0, y_1,\ldots,y_n$
are contained in $S_{\varphi}[u]$. Denote by $\<u\>$ the principal ideal in
$S_{\varphi}[u]$ generated by $u$. Using \cite[Proposition~1.2]{CDS}, we
obtain the formula  
\[
  H_f = \sum_{J\subset I} (V(J))^2 \prod_{j\in J} y_j,
\]
where the sum runs over all subsets $J = \{j_0,j_1,\ldots,j_d\}$ in $I =
\{0,1,\ldots,n\}$ and  $V(J)$ is the normalized $d$-dimensional volume of the
convex hull $conv(\{v_{j_0},v_{j_1},\ldots,v_{j_d}\})$ (in particular, 
$V(J) = 0$ if dimension of $conv(\{v_{j_0},v_{j_1},\ldots,v_{j_d}\})$ is less than
$d$). If $v_{i_1}, \ldots, v_{i_l}$ are not vertices of any $d$-dimensional
simplex $\tau_i\in {\mathcal T}$, then 
\begin{eqnarray}
\label{y}
  y_{i_1} \cdots  y_{i_l} = 
  u^{\varphi(v_{i_1}) + \cdots + \varphi(v_{i_l})} t_0^l 
  t^{v_{i_1}  + \cdots +         v_{i_l}} = 
\end{eqnarray}
\[
  = u^{\varphi(v_{i_1}) + \cdots + \varphi(v_{i_l})- 
  \varphi(v_{i_1}  + \cdots +         v_{i_l})} 
  \left(u^{\varphi(v_{i_1}  + \cdots +         v_{i_l})}t_0^l 
  t^{v_{i_1}  + \cdots +         v_{i_l}}
  \right) \in \<u\>,
\]
because the strict convexity of $\varphi$ implies that 
\[ 
  \varphi(v_{i_1}) + \cdots + \varphi(v_{i_l}) > 
  \varphi(v_{i_1}  + \cdots +         v_{i_l}).
\]
Therefore, the Hessian $H_f$ can be written as 
\[ 
  H_f = \sum_{i = 1}^k ({\rm Vol}(\tau_i))^2 
  \prod_{j : v_j \in \tau_i} y_j \quad + \quad h 
\]
for some  $h\in \< u \>$. Since every simplex $\tau_i\in {\mathcal T}$
contains $v_0 = 0$ and $H_f$ is divisible by $y_0$, then 
\begin{equation}
\label{H'}
  H_f' =  \sum_{i = 1}^k ({\rm Vol}(\tau_i))^2 
  \prod_{j : v_j \in \tau_i, j\neq 0} y_j \quad + \quad h'
\end{equation}
for some $h'\in \< u \>$.

Let $e_1,\ldots,e_d$ be any basis of the dual lattice $N = Hom(M,\Z)$. 
Then we can write every monomial $t^m$ ($m\in M$) as product 
$t_1^{m_1}\cdots t_d^{m_d}$, where $m_i := \<m,e_i\>$. Denote by $\C(u)$ the
field of rational functions in variable $u$. The toric residue over $\C(u)$ is
uniquely determined by $\C[u]$-linear mapping 
\[
  {\rm Res}_f[u]: S_\Delta^d[u] \rightarrow \C(u),
\]
having following two properties: 

$(1)$ ${\rm Res}_f[u](H_f') = {\rm Vol}(\Delta)$; 

$(2)$ ${\rm Res}_f[u]$ vanishes on all $\C[u]$-submodules 
$F_i S_\Delta^{d -  1}[u]\subset S_\Delta^d[u]$ 
($0\le i \le d$), where 
\[
  F_0 := t_0 f(t) = y_0 + \sum_{l = 1}^n y_l, 
\]
\[
  F_i := t_0 t_i \partial f /\partial t_i = 
\sum_{l = 1}^n \<v_l,e_i\> y_l, \quad 
  i = 1, \ldots, d.
\]

Since the toric residue ${\rm Res}_f$ is the specialization of the
$\C[u]$-linear mapping ${\rm Res}_f[u]$ at the point $u = 0$, 
$R_P(u)= {\rm Res}_f[u](P(y_1,\ldots,y_n)) \in \Q(u)$ is regular at the point
$u = 0$ and 
\[
  R_P(0) = (-1)^d \lim_{u \to 0} {\rm Res}_f(t_0^d P(u^{\varphi(v_1)}t^{v_1}, 
  \ldots,u^{\varphi(v_n)}t^{v_n})
\]
for any homogeneous polynomial $P(y_1,\ldots,y_n)\in \Q[y_1,\ldots,y_n]$ of
degree $d$. 

The cohomology ring $H^*(\P,\Q)$ of the projective simplicial toric variety
$\P$ can be computed as a quotient of the polynomial ring $\Q[y_1,\ldots,y_n]$
by the sum of two ideals: $\<F_1,\ldots,F_d\>$ and the ideal generated by all
monomials  $y_{i_1} \cdots  y_{i_l}$ such that $v_{i_1}, \ldots, v_{i_l}$ are
not vertices of any $d$-dimensional simplex $\tau_i\in {\mathcal T}$. In this 
description,  the variables $y_1,\ldots,y_n$ represent  the classes of
Weil divisors $D_1,\ldots,D_n$ in $H^2(\P,\Q)$. Moreover, if $v_{i_1}, \ldots,
v_{i_d}$ are the vertices of a $d$-dimensional simplex $\tau_i\in {\mathcal T}$,
then the intersection number $[D_{i_1}]\cdots [D_{i_d}]$ equals 
$1/{\rm Vol}(\tau_i)$. 

It follows from (\ref{y}) and from the property $(2)$ of ${\rm Res}_f[u]$ that
two linear maps 
\[
  P(y_1,\ldots,y_n)\mapsto \int_\P P([D_1],\ldots,[D_n])
\]
and
\[
  P(y_1,\ldots,y_n)\mapsto R_P(0)
\]
have the same kernel in the space of homogeneous polynomials of degree
$d$. Therefore, in order to identify these linear maps, it is sufficient to
compare their values on the special polynomial 
\[
  \widetilde{P}(y_1,\ldots,y_n) := \sum_{i = 1}^k ({\rm Vol}(\tau_i))^2 
  \prod_{j : v_j \in \tau_i, j\neq 0} y_j.
\]
By the property $(1)$ of ${\rm Res}_f[u]$ and (\ref{H'}), we obtain 
$R_{\widetilde{P}}(0) = {\rm Vol}(\Delta)$. On the other hand, 
\begin{eqnarray*}
  \int_{\P} \widetilde{P}([D_1],\ldots,[D_n]) =  
  \sum_{i = 1}^k ({\rm Vol}(\tau_i))^2 \int_{\P} 
  \prod_{j : v_j \in \tau_i, j\neq 0} [D_j]   \\ = 
  \sum_{i = 1}^k ({\rm Vol}(\tau_i))^2 \frac{1}{{\rm Vol}(\tau_i)} = 
  \sum_{i = 1}^k {\rm Vol}(\tau_i) = {\rm Vol}(\Delta).
\end{eqnarray*}
This finishes the proof. 
\end{proof}

\begin{rem} We remark that the number ${\rm Vol}(\Delta)$ equals the stringy
  Euler number $e_{\rm st}(\P)$ of the simplicial toric variety $\P$ 
  \cite{Batyrev3}. It is known that the usual Euler number $e(\P)$
  equals $k$ (the number of $d$-dimensional cones in $\Sigma$). The top Chern
  class of $\P$   is represented in the cohomology ring $H^*(\P, \Q)$ by the polynomial 
 \[ C_d(y_1, \ldots, y_n) := \sum_{i = 1}^k {\rm Vol}(\tau_i) 
  \prod_{j : v_j \in \tau_i, j\neq 0} y_j. \]
In the proof of the last theorem, we have shown that the
Hessian $H_f'$ specializes at $u =0$ to the polynomial  
$\widetilde{P}(y_1,\ldots,y_n)$ which represents the stringy top Chern
 class of $\P$.  
\end{rem}
\bigskip

\section{Hirzebruch surface ${\mathbb F}_1$}
\label{Sect_Hirz}

We consider below 
 a simplest example which illustrates many interesting 
ingredients of our conjecture.

\begin{figure}[ht]
  \centering \includegraphics*[scale=0.5]{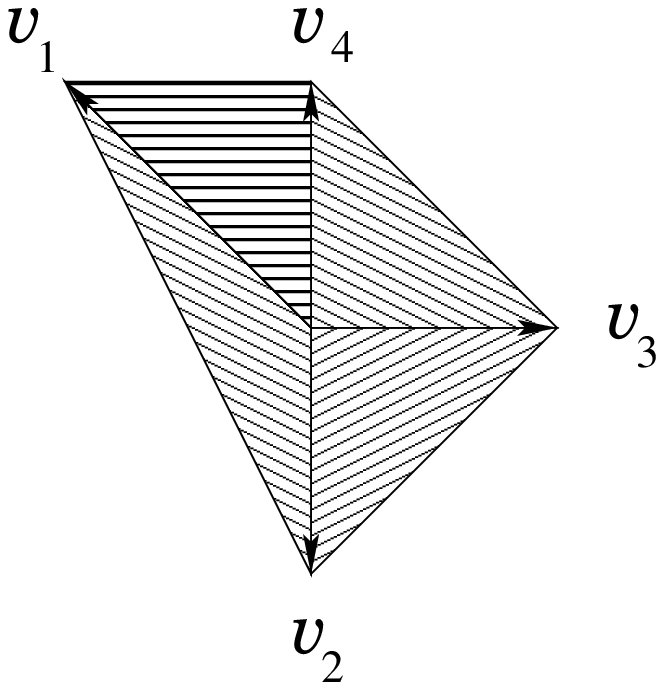}
  \caption{Fan ${\Si}$ for the Hirzebruch surface $\P_\Si = {\mathbb F}_1$}
  \label{Fan_Hirz}
\end{figure}

Let $M\cong \Z^2$ and $\Delta\subset M_\R\cong\R^2$ be a reflexive polytope
with the vertices 
\[
  v_1 = (-1,1),\quad v_2 = (0,-1),\quad v_3 = (1,0),\quad v_4 = (0,1).
\]
Denote by $A$ the set of points $\{0,v_1,\ldots,v_4\}$. Take a coherent
triangulation ${\mathcal T} = \{\tau_1,\ldots,\tau_4\}$ of $\Delta$ associated
with $A$ such that $0$ is a vertex of every $\tau_i$. The surface 
${\mathbb F}_1$ is a toric variety $\P_\Si$ defined by  the fan 
$\Si = \Si({\mathcal T})\subset M_\R$ (see Figure~\ref{Fan_Hirz}) whose
$1$-dimensional cones are generated by $v_i$. 
Take a generic Laurent polynomial 
\[
  f(t) = 1 - a_1 t^{v_1} - a_2 t^{v_2} - a_3 t^{v_3} - a_4 t^{v_4} =
         1 - a_1 t_1^{-1}t_2 - a_2 t_2^{-1} - a_3 t_1 - a_4 t_2
\]
with support in $\Delta$. 

\begin{prop}
\label{Princ_E_A_Hirz}
The principal $A$-determinant of $f(t)$ is 
\begin{eqnarray*}
  E_A(f) = 
     \underline{a_1^2 a_2^2 a_3^2 a_4^2}
   + \underline{a_1^3 a_2^2 a_3^3 a_4}   
            -8  a_1^2 a_2^3 a_3^2 a_4^3 
 +\underline{16 a_1^2 a_2^4 a_3^2 a_4^4}
            -36 a_1^3 a_2^3 a_3^3 a_4^2 
 -\underline{27 a_1^4 a_2^3 a_3^4 a_4}. 
\end{eqnarray*}
The underlined terms are in one-to-one correspondence with the vertices of the 
secondary polytope ${\rm Sec}(A)$. Moreover, the vertex $v_{\mathcal T}$
corresponding to the triangulation ${\mathcal T}$ is related with the monomial
$a_1^2 a_2^2 a_3^2 a_4^2$. 
\end{prop}

\begin{figure}[ht]
  \centering \includegraphics*[scale=0.5]{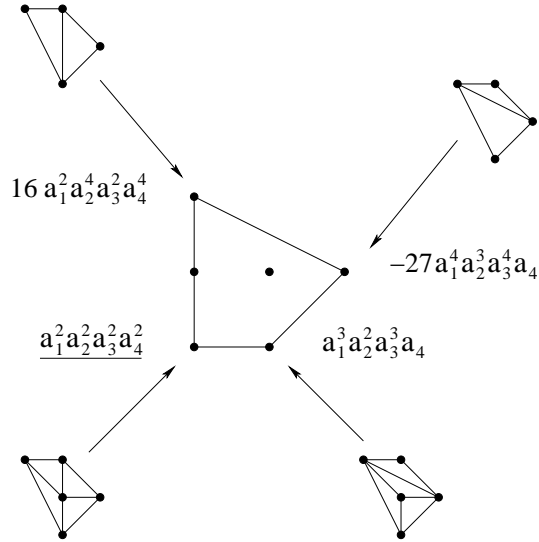}
  \caption{Secondary polytope ${\rm Sec}(A)$}
  \label{Sec_F_1}
\end{figure}

\begin{proof} 
It is easy to find all coherent triangulations of $\Delta$ (see 
Figure~\ref{Sec_F_1}). By knowing the coherent triangulations, using
Theorem~\ref{E_A=Sec(A)}, we can compute (up to sign) the terms of
$E_A(f)$ corresponding to the vertices of ${\rm Sec}(A)$: 
\[
  a_1^2 a_2^2 a_3^2 a_4^2, \,
  a_1^3 a_2^2 a_3^3 a_4, \,
  16 a_1^2 a_2^4 a_3^2 a_4^4, \,
  27 a_1^4 a_2^3 a_3^4 a_4.
\] 

In order to find the other terms, we compute (up to multiplication by monomial)
 the discriminant $D_A(f)$,
i.e., the set of those coefficients $\{a_i\}$ of $f(t)$ such that the system 
\begin{equation*}
\label{sys_Hirz}
  f(t) = t_1 \frac{\partial f}{\partial t_1}(t) = 
         t_2 \frac{\partial f}{\partial t_2}(t) = 0
\end{equation*}
has a solution in the torus $\T \cong (\C^*)^2$. If we define  
\[
  Z_1 := t_1^{-1} t_2, \, Z_2 := t_2^{-1}, \, Z_3 := t_1, \, Z_4 := t_2,
\]
then we obtain $Z_4 = Z_1 Z_3$, $Z_2 Z_4 = 1$ and previous three equations 
can be rewritten as 
\begin{eqnarray*}
  1 - a_1 Z_1 - \cdots - a_4 Z_4 = 0, \, 
  a_1 Z_1 - a_3 Z_3 = 0,\\ a_2 Z_2 - a_1 Z_1 - a_4 Z_4 = 0.
 \end{eqnarray*}
Excluding $Z_3$ and $Z_4$, we get 
\[
  1 - a_1 Z_1 - 2a_2 Z_2  = 0, \,
  a_3 (a_2 Z_2 - a_1 Z_1) = a_1 a_4 Z_1^2,\,
  Z_2 (a_2 Z_2 - a_1 Z_1) = a_4.
\]
The last system is equivalent to two  homogeneous equations 
\[
  a_3 (a_1 Z_1 + 2a_2 Z_2)(a_2 Z_2 - a_1 Z_1) = a_1 a_4 Z_1^2,\,
  Z_2 (a_2 Z_2 - a_1 Z_1) = a_4(a_1 Z_1 + 2a_2 Z_2)^2.
\]
Computing the resultant of two polynomials, we get the discriminant
\[
  D_A(f) = a_4
         + a_1 a_3  
      - 8  a_2 a_4^2 
      + 16 a_2^2 a_4^3
      - 36 a_1 a_2 a_3 a_4 
      - 27 a_1^2 a_2 a_3^2.
\]
We remark that 
the Newton polytope of $E_A(f)$ (secondary polytope) is obtained 
from the Newton polytope of $D_A(f)$ via a shift by 
vector $(2,2,2,1)$. So the multiplication of $D_A(f)$ by 
$a_1^2 a_2^2 a_3^2 a_4$ 
yields a polynomial having ${\rm Sec}(A)$ as its  Newton polytope
${\rm Sec}(A)$. Since $D_A(f)$ divides $E_A(f)$ (see 
Remark~\ref{discrim}), we conclude that $a_1^2 a_2^2 a_3^2 a_4 D_A(f)$ 
is exactly $E_A(f)$. 
\end{proof}

For the case $\P = {\mathbb F_1}$,
Conjecture~\ref{TRMC} can be reformulated as follows: 

\begin{conj} 
\label{Conj_F_1}
Let $D_1,\ldots,D_4$ be the torus-invariant divisors on $\P_\Si$ corresponding
to the vectors $v_1,\ldots,v_4$ respectively. Fix a homogeneous polynomial 
$P(x_1,\ldots,x_4)$ of degree two in $\Q[x_1,\ldots,x_4]$. Then the series
expansion of the toric residue 
\[ 
  R_P(a_1,\ldots,a_4) = 
  {\rm Res}_f (t_0^2\, P(a_1 t^{v_1}, \ldots,a_4 t^{v_1}))
\]
at the vertex $v_{\mathcal T}\in {\rm Sec}(A)$ defined by the triangulation
${\mathcal T}$ coincides with the generating function of the 
intersection numbers 
\[
  I_P(a_1,\ldots,a_4) = 
  \sum_{\beta\in K_{\rm eff}(\P_\Si)} \<P([D_1],\ldots,[D_4])\, \Phi_\beta\>_\beta
  \,a^\beta,
\]
where
$ \Phi_\beta = ([D_1] + \cdots + [D_4])^{b_1 + \cdots + b_4} \prod_{i : b_i < 0}
[D_i]^{-b_i - 1}$, $a^\beta = a_1^{b_1}\cdots a_4^{b_4}$. 
\end{conj}

It is possible to compute explicitly both the generating function for
intersection numbers and the toric residue. In this way, we check the equality
of Conjecture~\ref{Conj_F_1} by direct calculation. We will omit the details of
this calculation and sketch only ideas of how it can be done. 

First, we remark that the Mori cone $K_{\rm eff}(\P_\Si) \subset 
R(\Sigma)_\R$ 
is spanned by two elements 
\[
  l^{(1)} = (1,0,1,-1),\quad l^{(2)} = (0,1,0,1),
\]
i.e., $\beta$ runs over all lattice points 
\[
  (b_1,\ldots,b_4) = 
  \lambda_1 (1,0,1,-1) + \lambda_2 (0,1,0,1) = 
  (\lambda_1,\lambda_2,\lambda_1,\lambda_2 - \lambda_1),\quad 
  \lambda_1,\lambda_2\ge 0.
\]
The dual to $K_{\rm eff}(\P_\Si)$ K\"ahler cone 
is generated by the classes of $D_1$ and $D_2$. 
There are two independent linear relations
between the classes of torus-invariant divisors $D_1, D_2,D_3,D_4$: 
\[
  [D_1] - [D_3] = 0, \quad [D_1] - [D_2] + [D_4] = 0.
\]
We can consider the classes of $D_1$, $D_2$ as 
generators of $H^2(\P_\Si, \Z)$ and put  
\[
  y_1 := a^{l^{(1)}} = \frac{a_1 a_3}{a_4}, \quad 
  y_2 := a^{l^{(2)}} = a_2 a_4
\]
so that $a^\beta = a_1^{b_1}\cdots a_4^{b_4}$ is equal to 
$y_1^{\lambda_1}y_2^{\lambda_2}$. Substituting  above linear 
relations, we get 
\[
  P([D_1],\ldots,[D_4]) = P([D_1],[D_2],[D_1],[D_2] - [D_1]).
\]
By linearity, it is sufficient to compute the generating functions
\[
  I_{x_1^i x_2^{2 - i}}(y_1,y_2) = 
  \sum_{\lambda_1,\lambda_2\ge 0} \<[D_1]^i[D_2]^{2 - i}\, \Phi_\beta\>_\beta
  \,y_1^{\lambda_1} y_2^{\lambda_2},\quad i = 0,1,2.
\]

\begin{figure}[ht]
  \centering \includegraphics*[scale=0.75]{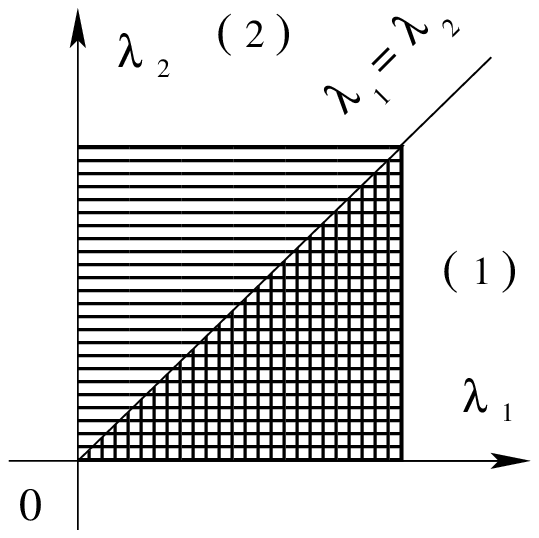}
  \caption{Two parts of the Mori cone of $\P_\Si$}
  \label{parts_Mori}
\end{figure}

We divide the Mori cone of $\P_\Si$ into two parts (see
Figure~\ref{parts_Mori}): 

\smallskip

$(1)$ $\lambda_1 > \lambda_2 \ge 0$; 

\smallskip

$(2)$ $\lambda_2 \ge \lambda_1 \ge 0$. 

\smallskip

For each of the part of the Mori cone we find the coefficients of the
generating function $I_{x_1^i x_2^{2 - i}}(y_1,y_2)$. 

\bigskip

{\sc Case 1.} 
If $\lambda_1 > \lambda_2 \ge 0$, then the last coordinate $b_4$ 
of 
\[
  \beta = (b_1,\ldots,b_4) = 
  (\lambda_1,\lambda_2,\lambda_1,\lambda_2 - \lambda_1)
\] 
is negative. By definition of the 
Morrison-Plesser moduli spaces $\P_\beta$ and the 
Morrison-Plesser classes $\Phi_\beta$, we obtain 
$\P_\beta = \P^{2\lambda_1 + 1}\times\P^{\lambda_2}$ and
\[
  \Phi_\beta = 
  ([D_1] + 2[D_2])^{\lambda_1 + 2\lambda_2}
  ([D_2] - [D_1])^{\lambda_1 - \lambda_2 - 1}.
\]
The intersection theory on $\P_\beta$ shows immediately that 
\[
  \<[D_1]^{l_1} [D_2]^{l_2}\>_\beta = 
  \left\{
  \begin{array}{ll}
  1,   & l_1 = 2\lambda_1 + 1, \quad l_2 = \lambda_2, \\
  0,   & {\rm otherwise},
  \end{array}
  \right.
\]
So we obtain 

\begin{align*}
  \<[D_1]^i[D_2]^{2 - i}\, \Phi_\beta\>_\beta = 
  \<[D_1]^i[D_2]^{2 - i} 
  ([D_1] + 2[D_2])^{\lambda_1 + 2\lambda_2}
  ([D_2] - [D_1])^{\lambda_1 - \lambda_2 - 1} 
  \>_\beta \\ = 
  \sum_{k \ge 0}
  (-1)^k 2^{2\lambda_2 - \lambda_1 + k - 1 + i} 
  \binom{\lambda_1 - \lambda_2 - 1}{k} 
  \binom{\lambda_1 + 2\lambda_2}{2\lambda_1 - k + 1 - i}. 
\end{align*}

\begin{rem} 
We remark that there is no irreducible curve $C
\subset {\mathbb F}_1$ such that $[C] = \beta = (b_1,\ldots,b_4)$, if 
$b_4 < 0$, and $b_2> 0$. Indeed, if $C \subset {\mathbb F}_1$ were such a
curve, then $b_4 = (C,D_4)<0$ would imply that $[C]$ is proportional to
$[D_4]$ (the latter  contradicts  $b_2 = (C,D_2) >0$ because of 
$(D_4,D_2) = 0$). On the other hand, the corresponding Morrison-Plesser
moduli spaces $\P_{\beta} = \P^{2b_1 + 1}\times\P^{b_2}$ are always nonempty.
\label{irred-cl}
\end{rem}

\bigskip

{\sc Case 2.} 
Let $\lambda_2 \ge \lambda_1 \ge 0$. Then the Morrison-Plesser moduli spaces 
$\P_\beta$ are toric varieties of
dimension  $\lambda_1 + 2\lambda_2 + 2$. The Morrison-Plesser class is equal
to $\Phi_\beta = ([D_1] + 2 [D_2])^{\lambda_1 + 2\lambda_2}$ and the
coefficients of the series $I_{x_1^i x_2^{2 - i}}(y_1,y_2)$ are the
intersection numbers 
\[
  \<[D_1]^i[D_2]^{2 - i}\,\Phi_\beta\>_\beta = 
  \<[D_1]^i[D_2]^{2 - i}([D_1] + 2[D_2])^{\lambda_1 + 2\lambda_2}\>_\beta. 
\]
These numbers can be computed directly using the intersection theory on the
Morrison-Plesser moduli spaces $\P_\beta$. We omit the details of the
proof. It is remarkable that the obtained formula 
\[
  \<[D_1]^i[D_2]^{2 - i}\,\Phi_\beta\>_\beta = 
  \sum_{k \ge 0}
  2^{2\lambda_2 - \lambda_1 + k - 1 + i} 
  \binom{\lambda_2 - \lambda_1 + k}{k} 
  \binom{\lambda_1 + 2\lambda_2}{2\lambda_1 - k + 1 - i}
\]
can be found from the formula obtained in Case 1
by the analytic continuation of the binomial coefficients: 
\[
  \binom{-m}{n} := (-1)^n\binom{m + n - 1}{n},\quad n = 0, \pm 1, \pm 2, \ldots
\]
for any positive integer $m$. 

\bigskip 

Summarizing all the considered cases, we obtain 
the generating function of
the intersection numbers in the form 
\begin{eqnarray*}
  &&I_{x_1^i x_2^{2 - i}}(y_1,y_2) = \\
  = 
  &&\sum_{\lambda_1, \lambda_2 \ge 0} 
  \left(\sum_{k \ge 0}
  (-1)^k 2^{2\lambda_2 - \lambda_1 + k - 1 + i} 
  \binom{\lambda_1 - \lambda_2 - 1}{k} 
  \binom{\lambda_1 + 2\lambda_2}{2\lambda_1 - k + 1 - i} \right)
  \,y_1^{\lambda_1} y_2^{\lambda_2}. 
\end{eqnarray*}

In order to compare this series with toric residues, we use the following 
integral representation: 
\[
  I_{x_1^i x_2^{2 - i}}(y_1,y_2) = 
  \sum_{\lambda_1, \lambda_2 \ge 0} 
  \left(
  \frac{1}{(2\pi i)^2}\int_\gamma
  \frac{u_0^i u_1^{2 - i}(u_1 - u_0)^{\lambda_1 - \lambda_2 - 1} 
  (u_0 + 2u_1)^{\lambda_1 + \lambda_2} du_0\wedge du_1}
  {u_0^{2\lambda_1 + 2} u_1^{\lambda_2 + 1}}
  \right)\,y_1^{\lambda_1} y_2^{\lambda_2},
\]
where $\gamma = \{(u_0,u_1)\in \C^2: |u_0| = \varepsilon_1, 
|u_1| = \varepsilon_2\}$, $\varepsilon_1, \varepsilon_2 > 0$. 
Changing the order of summation and integration, we get 
\[
  I_{x_1^i x_2^{2 - i}}(y_1,y_2) = 
  \frac{1}{(2\pi i)^2}
  \int_\gamma
  \frac{u_0^i u_1^{2 - i} du_0\wedge du_1}
  {(u_0^2 - (u_0 + 2u_1)(u_1 - u_0)y_1)(u_1(u_1 - u_0) - (u_0 + 2u_1)^2 y_2)}.
\]
This integral can be computed directly. We write down the obtained rational
functions $I_Q(y_1,y_2)$ 
together with some first terms of their expansions:

\begin{itemize}

\smallskip

\item $Q(x_1,x_2) = x_2^2$:
\begin{multline*}
  \frac{1 + y_1 + 4y_2 + 3y_1 y_2}
  {1 + y_1 - 8 y_2 + 16 y_2^2 - 36 y_1 
  y_2 - 27 y_1^2 y_2} = 
  1 + 12 y_2 + 27 y_1 y_2 + 80 y_2^2 \\ 
  + 568 y_1 y_2^2 + 448 y_2^3 + 728 y_1^2 y_2^2 + 
  6544 y_1 y_2^3 + 2304 y_2^4 + y_1^3 y_2^2 + 21888 y_1^2 y_2^3 + \cdots
\end{multline*}

\item $Q(x_1,x_2) = x_1x_2$ :
\begin{multline*}
  \frac{1 + y_1 - 4y_2 - 6y_1 y_2}
  {1 + y_1 - 8 y_2 + 16 y_2^2 - 36 y_1 
  y_2 - 27 y_1^2 y_2} = 
  1 + 4 y_2 + 26 y_1 y_2 + 16 y_2^2 
  + y_1^2 y_2 \\ + 336 y_1 y_2^2 
  + 64 y_2^3 - y_1^3 y_2 + 716 y_1^2 y_2^2 + 2784 y_1 y_2^3 + 
  256 y_2^4 + y_1^4 y_2 + 14 y_1^3 y_2^2 + \cdots
\end{multline*}

\item  $Q(x_1,x_2) = x_1^2$: 
\begin{multline*}
  \frac{y_1(1 + 12 y_2)}
  {1 + y_1 - 8 y_2 + 16 y_2^2 - 36 y_1 
  y_2 - 27 y_1^2 y_2} = 
  y_1 - y_1^2 + 20 y_1 y_2 +  
  y_1^3 + 8 y_1^2 y_2  \\ + 144 y_1 y_2^2 
  - y_1^4 - 9 y_1^3 y_2 + 656 y_1^2 y_2^2 + 
  832 y_1 y_2^3 + y_1^5 + 10 y_1^4 y_2 +
  84 y_1^3 y_2^2 + \cdots
\end{multline*}

\end{itemize}

These rational functions can be identified with the toric residues 
$R_{x_1^ix_2^{2 - i}}(y_1,y_2)$ $(i=0,1,2)$ which one computes, for example,
by the method from Section~\ref{Sect_Yukawa_CY}. This verifies
Conjecture~\ref{Conj_F_1} directly.

\bigskip

\section{Toric residue and flop}

In this section, we consider a  simplest reflexive polytope $\Delta$ which has
two different coherent triangulations such that $0$ is a vertex of all
simplices. These triangulations correspond to two different vertices of the
secondary polytope. We compute the corresponding expansions of the toric
residues at each of these vertices. 

Denote by $A$ the union of the origin $v_0 = (0,0,0)$ in $M\cong\Z^3$ together
with the points 
\[
  v_1 = (1,0,0),\, v_2 = (0,1,0),\, v_3 = (-1,-1,0),\, 
  v_4 = (0,0,1),\, v_5 = (1,1,-1).
\]
Then $\Delta := conv(A)\subset M_{\R}\cong\R^3$ is the reflexive polytope. 
Let 
\[
  f(t) = 1 - \sum_{i = 1}^5 a_i t^{v_i} = 
  1 - a_1 t_1 - a_2 t_2 - a_3 t_1^{-1} t_2^{-1} - a_4 t_3 - 
  a_5 t_1 t_2 t_3^{-1}
\]
be a generic Laurent polynomial.

\begin{prop}
The principal $A$-determinant of $f(t)$ equals 
\begin{eqnarray*}
   E_A(f) = 
  \underline{a_1^4 a_2^4 a_3^4 a_4^3 a_5^3} 
  - 54   a_1^5 a_2^5 a_3^5 a_4^3 a_5^3
  -\underline{a_1^3 a_2^3 a_3^4 a_4^4 a_5^4}
  +\underline{729  a_1^6 a_2^6 a_3^6 a_4^3 a_5^3}
   \\
  +54   a_1^3 a_2^3 a_3^5 a_4^5 a_5^5
  - 2187 a_1^5 a_2^5 a_3^6 a_4^4 a_5^4
  + 2187 a_1^4 a_2^4 a_3^6 a_4^5 a_5^5
  -\underline{729  a_1^3 a_2^3 a_3^6 a_4^6 a_5^6},
\end{eqnarray*}
where the terms corresponding to the vertices of the secondary polytope 
${\rm Sec}(A)$ are underlined. 
\end{prop}

\begin{proof} The idea of the proof is the same as in
Proposition~\ref{Princ_E_A_Hirz}. The terms corresponding to the vertices of
the polytope ${\rm Sec}(A)$ can be easily found (up to sign) from the 
coherent triangulations of the polytope $\Delta$: 
\[
      a_1^4 a_2^4 a_3^4 a_4^3 a_5^3,\,
      a_1^3 a_2^3 a_3^4 a_4^4 a_5^4,\,
  729 a_1^6 a_2^6 a_3^6 a_4^3 a_5^3,\,
  729  a_1^3 a_2^3 a_3^6 a_4^6 a_5^6.
\]

The Laurent polynomial  $f$ is  $\Delta$-regular, 
i.e., $E_A(f) \neq 0$, if and only if 
the equations 
\begin{equation*}
  f(t) = t_1 \frac{\partial f}{\partial t_1}(t) = 
         t_2 \frac{\partial f}{\partial t_2}(t) = 
         t_3 \frac{\partial f}{\partial t_3}(t) = 0
\end{equation*}
have no solution in the compactification  $\P_\Delta$ of the 
torus $\T\cong (\C^*)^3$. If we put 
\[
  Z_1 := t_1,\, Z_2 := t_2,\, Z_3 := t_1^{-1} t_2^{-1},\, Z_4 := t_3,\,
  Z_5 := t_1 t_2 t_3^{-1},
\]
then the $Z_0$-homogenization of the last system is equivalent to 
\begin{eqnarray*}
  Z_0 - a_1 Z_1 - \cdots - a_5 Z_5 = a_4 Z_4 - a_5 Z_5 
  = a_1 Z_1 - a_3 Z_3 + a_5 Z_5 \\ 
  = a_2 Z_2 - a_3 Z_3 + a_5 Z_5 = 0;\,
  Z_1 Z_2 Z_3 = Z_0^3,\, Z_1 Z_2 = Z_4 Z_5.
\end{eqnarray*}

By excluding $Z_1, Z_2, Z_4$, we obtain 
\[
  Z_0 - 3 a_3 Z_3 = (a_3 Z_3 - a_5 Z_5)^2 Z_3 - a_1 a_2 Z_0^3 = 
  a_4 (a_3 Z_3 - a_5 Z_5)^2 - a_1 a_2 a_5 Z_5^2 = 0. 
\]
Hence, $\Delta$-regularity of $f$ is equivalent 
to nonvanishing of the resultant $R(g_1,g_2)$ of two homogeneous
forms 
\[
  g_1 = (a_3 Z_3 - a_5 Z_5)^2 Z_3 - 27 a_1 a_2 a_3^3 Z_3^3, \; 
g_2  = 
  a_4 (a_3 Z_3 - a_5 Z_5)^2 - a_1 a_2 a_5 Z_5^2.
\]
By   direct computation, we obtain that $R(g_1,g_2)$ splits into product 
of two irreducible divisors of multiplicity $1$. Moreover, it is 
straightforward to see that the
difference between the Newton polytopes 
of $R(g_1,g_2)$ and $E_A(f)$ is the vector 
$(1,1,0,3,0)$. Since $R(g_1, g_2)$ and $E_A(f)$ have the same set 
of irreducible divisors $D_{A\cap \Gamma}(f)$, $\Gamma \subseteq \Delta$   
(see Remark~\ref{discrim}), we see that  
the multiplication of $R(g_1,g_2)$ by 
$-a_1 a_2 a_4^3$ yields  $E_A(f)$. 
\end{proof}

The toric residue in the following statement can be computed using, for
example, the method from Section~\ref{Sect_Yukawa_CY}. 

\begin{prop}
Let $P(x_1,x_2,x_3,x_4) = x_1 x_2 x_4\in\Q[x_1,\ldots,x_4]$ be the input
polynomial. Then the corresponding toric residue 
\[
  R_{x_1 x_2 x_4}(a_1,\ldots,a_5) = 
  -{\rm Res}_f(t_0^3\,(a_1 t_1)(a_2 t_2)(a_4 t_3)) 
\]
can be expressed as a quotient 
\[
  R_{x_1 x_2 x_4}(a_1,\ldots,a_5) = 
   (    a_1^4 a_2^4 a_3^4 a_4^3 a_5^3 - 
     27 a_1^5 a_2^5 a_3^5 a_4^3 a_5^3 -
     81 a_1^4 a_2^4 a_3^5 a_4^4 a_5^4)/E_A(f).
\]
\end{prop}

\begin{figure}[h]
  \includegraphics*[width=0.37\textwidth]{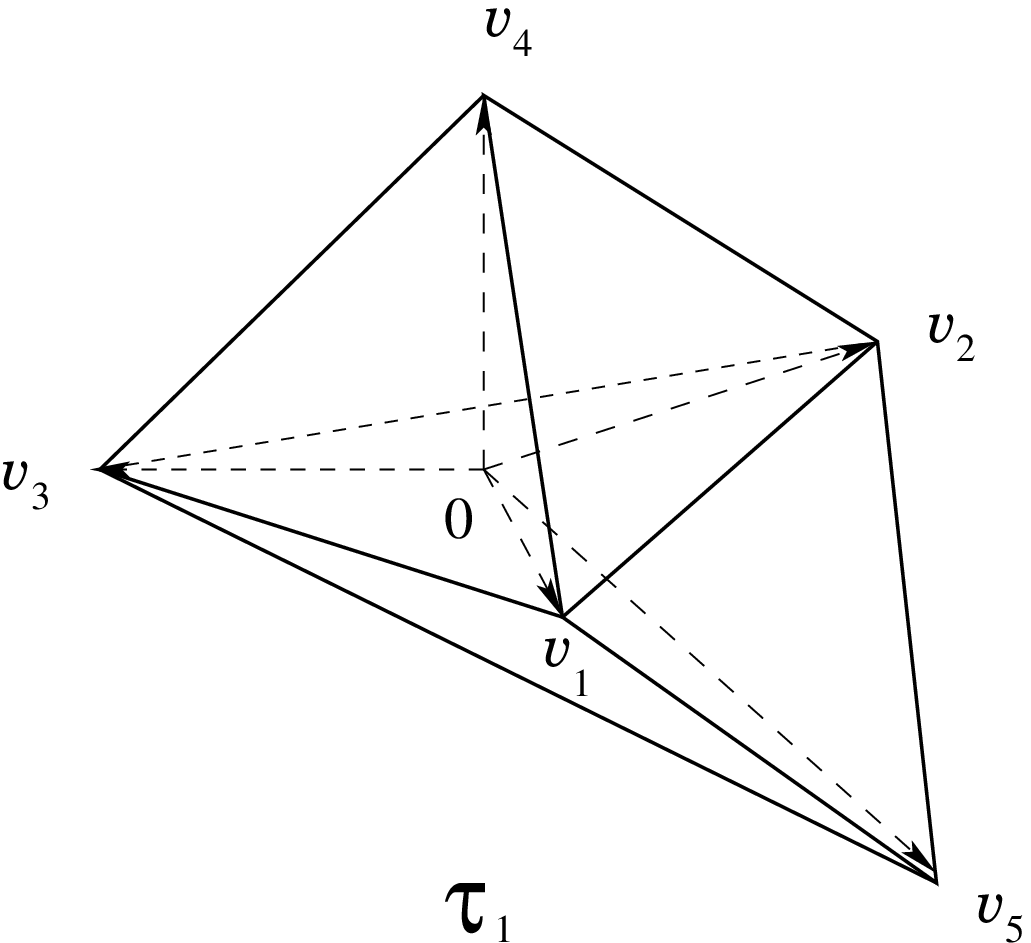}
  \hfill
  \includegraphics*[width=0.37\textwidth]{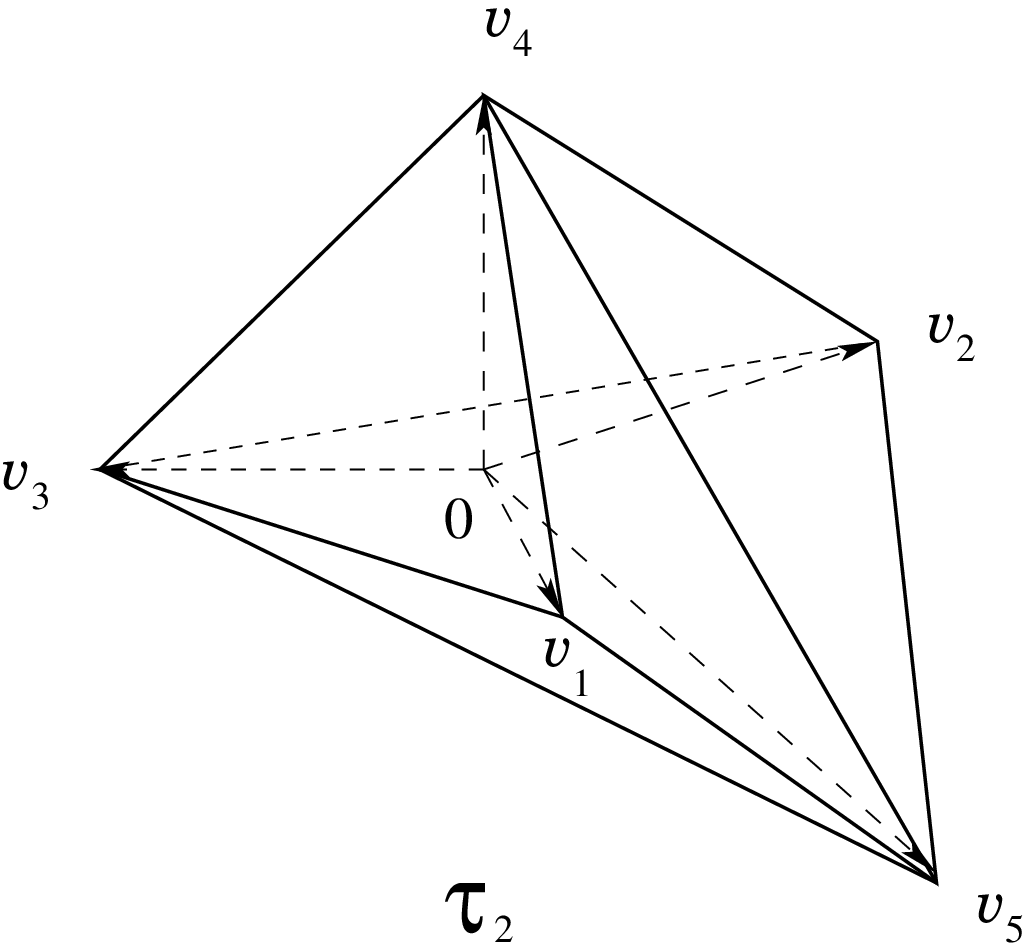}
  \\
  \caption{Triangulations ${\mathcal T}_1$ and ${\mathcal T}_2$} 
  \label{triang_Flop}
\end{figure}

There are two triangulations related with the vectors $v_1, \ldots, v_5$ with
the property that their maximal-dimensional simplices contain $0$ (see 
Figure~\ref{triang_Flop}). Namely, 
\[
  {\mathcal T}_1 = \{\tau_{124},\, \tau_{125},\, \tau_{134},\, \tau_{135},\, 
  \tau_{235},\, \tau_{234}\}
\]
and 
\[
  {\mathcal T}_2 = \{\tau_{145},\, \tau_{245},\, \tau_{134},\, \tau_{135},\, 
  \tau_{235},\, \tau_{234}\}, 
\]
here $\tau_{124}$ means that this simplex is generated by $0$, $v_1$, $v_2$,
$v_4$, etc. Let $\Si_1 = \Si({\mathcal T}_1) \subset M_\R$ 
($\Si_2 = \Si({\mathcal T}_2) \subset M_\R$) be the fan whose $d$-dimensional
simplices are defined as $\sigma := \R_{\ge 0}\tau$, 
$\tau\in {\mathcal T}_1$ ($\tau\in {\mathcal T}_2$). Toric varieties
$\P_{\Si_1}$ and $\P_{\Si_2}$ corresponding to the fans $\Si_1$ and $\Si_1$
are related by a {\it flop}. 

\begin{figure}[ht]
  \centering \includegraphics*[scale=0.6]{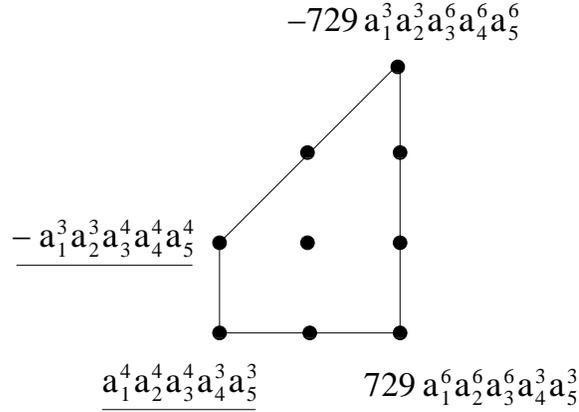}
  \caption{Secondary polytope ${\rm Sec}(A)$}
  \label{SecForFlop}
\end{figure}

The secondary polytope ${\rm Sec}(A)$ is depicted in Figure~\ref{SecForFlop}, 
where the monomials 
$a_1^4 a_2^4 a_3^4 a_4^3 a_5^3$ and $a_1^3 a_2^3 a_3^4 a_4^4 a_5^4$ of 
$E_A(f)$ corresponding to the vertices $v_{{\mathcal T}_1}$ and 
$v_{{\mathcal T}_2}$ of ${\rm Sec}(A)$ respectively are underlined. 

In oder to find the power series expansion of the toric residue 
$R_{x_1 x_2 x_4}(a_1,\ldots,a_5)$ at the vertex $v_{{\mathcal T}_1}$ of 
${\rm Sec}(A)$, we 
rewrite it  as follows 
\[
  R_{x_1 x_2 x_4}(u_1,u_2) = 
  \frac{1 - 27 u_1 - 81 u_1 u_2}
  {(1 - u_2)(1 - 54 u_1 + 729 u_1^2 - 54 u_1 u_2 - 
  1458 u_1^2 u_2 + 729 u_1^2 u_2^2)},
\]
where the variables   
\[
  u_1 := a_1 a_2 a_3,\, u_2 := \frac{a_4 a_5}{a_1 a_2}
\]
correspond to the generators 
\[
  l_1=(1,1,1,0,0), \;\; l_2=(-1,-1,0,1,1)
\]
of the Mori cone of $\P_{\Si_1}$. Here are the first terms of the Taylor
expansion
\begin{eqnarray*}
  R_{x_1 x_2 x_4}(u_1,u_2) = 
  1 + 27 u_1 + u_2 + 729 u_1^2 + u_2^2 + 19683 u_1^3 + 2187 u_1^2 u_2 + u_2^3
  + \cdots 
\end{eqnarray*}
 The power series expansion of the toric residue 
$R_{x_1 x_2 x_4}(a_1,\ldots,a_5)$ at the vertex $v_{{\mathcal T}_2}
 \in {\rm Sec}(A)$
(corresponding to the second triangulation ${\mathcal T}_2$)
needs another 
variables 
\[
  w_1 := u_1 u_2 = a_3 a_4 a_5,\, 
  w_2 := u_2^{-1} = \frac{a_1 a_2}{a_4 a_5},
\]
corresponding to the generators 
\[
  l_1'=(0,0,1,1,1), \;\; l_2'=(1,1,0,-1,-1)
\]
of the Mori cone of $\P_{\Si_2}$.
We rewrite $R_{x_1 x_2 x_4}(a_1,\ldots,a_5)$ in
the form 
\begin{eqnarray*}
  &&R_{x_1 x_2 x_4}(w_1,w_2) = \\ &&=
  \frac{w_2 - 27w_1 w_2^2 - 81 w_1 w_2}
  {(w_2 - 1)(1 - 54w_1 - 54w_1w_2 + 
  729w_1^2 - 1458w_1^2w_2 + 729 w_1^2 w_2^2)}.
\end{eqnarray*}

We note that the power series expansion of $R_{x_1 x_2 x_4}(u_1,u_2)$ begins
with $1$, whereas the series expansion of $R_{x_1 x_2 x_4}(w_1,w_2)$ does not
contain a constant term. 
This fact agrees with different  intersection numbers of three divisors 
$D_1,D_2,D_4$  on $\P_{\Si_1}$ and
$\P_{\Si_2}$ corresponding to the vectors $v_1,v_2,v_4$:  
$\int_{\P_{\Si_1}} [D_1] [D_2] [D_4] = 1$ and 
$\int_{\P_{\Si_2}} [D_1] [D_2] [D_4] = 0$.

\bigskip

\section{Weighted projective spaces}

Consider  the $d$-dimensional weighted projective 
space $\P(w_1,\ldots,w_n)$, $n = d + 1$, where 
${\rm gcd}(w_1,\ldots,w_n) = 1$ and 
\[
  w_i | (w_1 + \cdots + w_n), \quad i = 1,\ldots,n.
\]
Let $\{v_1,\ldots,v_n\}\subset M\cong \Z^d$ be the vectors 
generating $M$ and satisfying the relation 
\[
  w_1 v_1 + \cdots + w_n v_n = 0.
\]
If we set  $A := \{v_0 = 0,v_1,\ldots,v_n \}$,
then the polytope $\Delta := conv(A) \subset M_\R$ is a reflexive
simplex. There are exactly  
two coherent
triangulations of $\Delta$: 
the triangulation ${\mathcal T}_1$ coinciding with
the whole polytope $\Delta$, and the triangulation 
${\mathcal T} = {\mathcal T}_2$ consisting of the union of the $d$-dimensional
simplices 
\[
  \tau_i = conv\{0, v_1, \ldots,\widehat{v_i},\ldots v_n\}, 
  \quad i = 1,\ldots,n, 
\]
where $\widehat{v_i}$ means that $v_i$ is omitted. Note that 
${\rm Vol}(\tau_i) = w_i$ and ${\rm Vol}(\Delta) = \sum_{i = 1}^n w_i$. 
The fan $\Sigma = \Sigma({\mathcal T})\subset M_\R$ defining
$\P(w_1,\ldots,w_n)$ has the generators $\{v_1,\ldots,v_n\}$. Let 
\[
  f(t) := 1 - \sum_{i = 1}^n a_i t^{v_i}
\]
be a generic Laurent polynomial. It is easy to find its principal
$A$-determinant from Theorem~\ref{E_A=Sec(A)}. 
\begin{prop} The secondary polytope ${\rm Sec}(A)$ is an interval and the
principal $A$-determinant of $f(t)$ is equal (up to sign) to 
\[
  E_A(f) = \prod_{i = 1}^n w_i^{w_i} 
  a_i^{w_1 + \cdots + \widehat{w}_i + \cdots + w_n} - 
  \left(\sum_{i = 1}^n w_i\right)^{\sum_{i = 1}^n w_i}
  (a_1\cdots a_n)^{\sum_{i = 1}^n w_i},
\]
where $\widehat{w_i}$ means that $w_i$ is omitted. The first summand in $
E_A(f)$ corresponds to the triangulation ${\mathcal T}$. 
\end{prop}

\begin{theo} 
\label{wp_MP}
Let $P(x_1,\ldots,x_n)\in\Q[x_1,\ldots,x_n]$ be an arbitrary homogeneous
polynomial of degree $d$. Denote $y := a_1^{w_1} \cdots a_n^{w_n}$. Then the
toric residue 
\[
  R_P(a) = 
  (-1)^d\,{\rm Res}_f(t_0^d\, P(a_1 t^{v_1},\ldots,a_n t^{v_n}))
\]
is equal to the rational function 
\[
  R_P(a_1^{w_1} \cdots a_n^{w_n}) = 
  R_P(y) = 
  \frac{\nu\cdot P(w_1,\ldots,w_n)}
  {1 - \mu\, y}, 
\]
where
\[
  \nu := \frac{1}{w_1\cdots w_n}, \quad
  \mu := \left(
  \frac{\displaystyle 
  \left(\sum_{i = 1}^n w_i\right)^{\sum_{i = 1}^n w_i}}
  {\displaystyle \prod_{i = 1}^n w_i^{w_i}} 
  \right).
\]
\end{theo}

\begin{proof}
By Theorem~\ref{Prop_tor_res}$(i)$, the toric residue $R_P$ can be computed as
a sum over the critical points $\xi$ of $f$:
 \[ 
  R_P = (-1)^d \sum_{\xi \in V_f}
  \frac{P(a_1\xi^{v_1}, \ldots,a_n\xi^{v_n} )}{f(\xi) {H}_f^0(\xi)}.  
\]
If $\xi$  is a critical point of $f$, i.e., a solution of the system of
equations 
\[
  t_1\frac{\partial f}{\partial t_1}(t) = \cdots = 
  t_d\frac{\partial f}{\partial t_d}(t) = 0, \quad t\in (\C^*)^d,
\]
then 
\[
  a_1 \frac{\xi^{v_1}}{w_1} = \cdots = a_n \frac{\xi^{v_n}}{w_n} = z
\]
and
\begin{equation}
  z^{w_1 + \cdots + w_n} = 
  \left(\frac{a_1}{w_1}\right)^{w_1} \cdots \left(\frac{a_n}{w_n}\right)^{w_n}.
\label{w_pr}
\end{equation}
These relations simplify our computations. For example, we may write 
at the critical points of $f$:
\[
  f(\xi) = 1 - \left(\sum_{i = 1}^n w_i\right) z,
\]
and 
\[
  P(a_1 \xi^{v_1},\ldots,a_n \xi^{v_n}) = 
  P(w_1,\ldots,w_n) z^d.
\]
The value of the polynomial  
\[
  H_f^0(t) = \det\left(\left(t_i\frac{\partial}{\partial t_i}\right)
  \left(t_j\frac{\partial}{\partial t_j}\right) f\right)_{1\le i,j \le d}  
\]  
at a critical point $\xi$ of $f$ equals 
\[
  H_f^0(\xi)  = (-1)^d w_1 \cdots w_n 
\left(\sum_{i = 1}^n w_i\right) z^d.
\]
Since the summation over the critical points is 
equivalent to the summation over the roots of the equation (\ref{w_pr}), we get 
\begin{eqnarray*}
  && R_P(y)  =
  \sum_{z^{w_1 + \cdots + w_n} = 
  \left(\frac{a_1}{w_1}\right)^{w_1} \cdots
  \left(\frac{a_n}{w_n}\right)^{w_n}}
  \frac{P(w_1,\ldots,w_n)}
  {w_1\cdots w_n(\sum_{i = 1}^n w_i) (1 - (\sum_{i = 1}^n w_i)z)} 
   \\ && =
  \frac{P(w_1,\ldots,w_n)}{w_1\ldots w_n} 
  \sum_{b \ge 0} 
  \left(\displaystyle \sum_{i = 1}^n w_i \right)^{(\sum_{i = 1}^n w_i)b}  
  \left(\frac{a_1^{w_1}\cdots a_n^{w_n}}{w_1^{w_1}\cdots w_n^{w_n}}\right)^b
  \\ && = \frac{\nu\cdot P(w_1,\ldots,w_n)}{1 - \mu\, y},
\end{eqnarray*}
as required.
\end{proof}

Now we show that the expansion of $R_P(y)$ at the vertex 
$v_{\mathcal T}\in {\rm Sec}(A)$ corresponding to the triangulation 
${\mathcal T}$ coincides with the generating function of intersection numbers:

\begin{theo}
\label{weght_2}
Let $P(x_1,\ldots,x_n)$ be any 
homogeneous polynomial in $\Q[x_1,\ldots,x_n]$
of degree $d$.
The generating function of intersection numbers on the Morrison-Plesser
moduli spaces has the form 
\[
  I_P(y) = \nu\cdot P(w_1,\ldots,w_n) \sum_{b \ge 0} \mu^by^b = 
  \frac{\nu\cdot P(w_1,\ldots,w_n)}{1 - \mu y}.
\]
\end{theo}

\begin{proof}
The Morrison-Plesser moduli spaces $\P_\beta$ for the weighted projective
space $\P = \P(w_1,\ldots,w_n)$ are also weighted projective spaces of
dimension $(\sum_{i = 1}^n w_i) b + d$. Let $D_1,\ldots,D_n$ be the divisors
corresponding to the vectors $v_1,\ldots,v_n$. We have the following relations
between the torus-invariant divisors on $\P_\beta$ modulo rational
equivalence: 
\[
  \frac{[D_1]}{w_1} = \cdots = \frac{[D_n]}{w_n} = [D_0].
\]
 Then the Mori cone $K_{\rm eff}(\P)$ consists of $\beta = (b)$,
$b\ge 0$. The Morrison-Plesser class is exactly 
\[
  \Phi_\beta = 
  \left(
  \left(
  \sum_{i = 1}^n w_i
  \right)
  [D_0]
  \right)^{(\sum_{i = 1}^n w_i) b}, 
\]
and the generating function of intersection numbers on $\P_\beta$ can be
written as 
\[
  I_P(y) = 
  \sum_{b \ge 0} 
  \<
  P([D_1],\ldots,[D_n]) 
  \left(
  \left(
  \sum_{i = 1}^n w_i
  \right)
  [D_0]
  \right)^{(\sum_{i = 1}^n w_i) b} 
  \>_\beta \,
  y^b.
\]
Using $\< [D_0]^{(\sum_{i = 1}^n w_i) b + d}\>_\beta = 
1/(w_1^{w_1 b + 1}\cdots w_n^{w_n b + 1})$, we obtain 
\begin{eqnarray*}
  I_P(y) &=& 
  P(w_1,\ldots,w_n) \sum_{b \ge 0} 
  \left(\sum_{i = 1}^n w_i\right)^{(\sum_{i = 1}^n w_i) b}
  \< [D_0]^{(\sum_{i = 1}^n w_i) b + d}\>_\beta\, y^b \\
  &=& 
  P(w_1,\ldots,w_n) \sum_{b \ge 0} 
  \left(\sum_{i = 1}^n w_i\right)^{(\sum_{i = 1}^n w_i) b}
  \frac{1}{w_1^{w_1 b + 1}\cdots w_n^{w_n b + 1}}\; y^b \\ 
  &=&
  \nu\cdot P(w_1,\ldots,w_n) \sum_{b \ge 0} \mu^by^b = 
  \frac{\nu\cdot P(w_1,\ldots,w_n)}{1 - \mu y}.
\end{eqnarray*}
\end{proof}

\bigskip

\section{Product of projective spaces}

In this section we check the Toric Residue Mirror Conjecture in the case 
$\P = \P^{d_1}\times\cdots\times\P^{d_r}$. 

For all $j = 1,\ldots,r$, 
we set $ n_j := d_j + 1$  and denote by $M_j$ the free abelian group of rank 
$d_j$ generated by the elements $v_{j1}, \ldots, v_{jn_j}$ satisfying the
linear relation: 
\[ 
  v_{j1} + \cdots + v_{j n_j} = 0. 
\]
Let $M:= M_1 \oplus \cdots \oplus M_r$. 
Consider a  Laurent polynomial
\begin{eqnarray*}
 f(t) = 1 - \sum_{i_1 = 1}^{n_1} a_{1 i_1} t^{v_{1 i_1}} - \cdots -
            \sum_{i_r = 1}^{n_r} a_{r i_r} t^{v_{r i_r}}
\end{eqnarray*}
 with support in the reflexive polytope 
\[ 
  \Delta = 
  conv(\{v_{11},\ldots,  v_{1 n_1}, \ldots,  v_{r1},\ldots,  v_{r n_r}\}) 
  \subset M_\R.
\]
The fan $\Sigma\subset M_\R$ consisting of cones over all faces 
of $\Delta$ determines the  toric variety 
$\P = \P^{d_1}\times\cdots\times\P^{d_r}$ of dimension 
$d := d_1 + \cdots + d_r$. Consider $n := n_1 + \cdots + n_r$ variables 
$x_{ji}$ ($1 \leq i \leq n_j, 1\leq j \leq r $).  
We will use the following notations:
\[
  x^k := 
  x_{11}^{k_{11}} \cdots x_{1 n_1}^{k_{1 n_1}}  \cdots
  x_{r1}^{k_{r1}} \cdots x_{r n_r}^{k_{r n_r}},\]
\[ 
  k_j := k_{j 1} + \cdots + k_{j n_j}, \quad
u_j:= n_j^{n_j} a_{j1}\cdots a_{j n_j},  \quad j = 1,\ldots,r.
\]
Let $K= \C(u_1,\ldots,u_r)$ be the field of rational functions
in $u_1, \ldots, u_r$. Let $z_j$ be a root of the equation 
 $z_j^{n_j} = u_j$ ($1\le j \le r$). 
We obtain a finite Galois extension $L= \C(z_1,\ldots,z_r)$ of $K$ 
of degree $[L:K] = n_1 \cdots n_r$. 
One has the algebraic trace map 
\[ {\rm tr}_{L/K}\; : \; L \to K, \]
which can be defined by the formula
\[ {\rm tr}_{L/K} (g) = \sum_{\begin{subarray}{l}
                    z_1^{n_1} = u_1 \\
                    \cdots             \\
                    z_r^{n_r} = u_r
                  \end{subarray}} g(z_1, \ldots, z_r),  \]
where  the sum runs  over all roots of the system of equations 
$z_j^{n_j} = u_j$ $(1\le j \le r)$. 

\begin{theo}
\label{trace_prod_proj}
Let $P(x) = x^k$ be any monomial of degree $d$ 
(i.e., $\sum_{j = 1}^r k_j = d$). Then the  toric residue $R_P$ 
corresponding to the monomial  $x^k$ 
is the rational function in $u_1,\ldots,u_r$: 
\[
R_{x^k}(u) = n_1^{d_1 - k_1 - 1}\cdots n_r^{d_r - k_r - 1}
 {\rm tr}_{L/K}\left(             
  \frac{1}{z_1^{d_1 - k_1}
    \cdots z_r^{d_r - k_r}(1 - z_1 - \cdots - z_r)} \right).
\] 
\end{theo}

\begin{proof}
The proof is based on the formula from Theorem~\ref{Prop_tor_res}$(i)$ 
which express the toric residue $R_P$ as the following 
sum over the critical points
$\xi$ of $f$:
\[ 
  R_P = (-1)^d \sum_{\xi \in V_f}
  \frac{P(a_{11}\xi^{v_{11}}, \ldots,a_{rn_r}\xi^{v_{rn_r}} )}{f(\xi) {H}_f^0(\xi)}.  
\]
If $\xi$  is a critical point of $f$, then the vanishing 
of partial derivatives of $f$ implies: 
\[
  a_{j1}\xi^{v_{j1}} = \cdots = a_{j n_j}\xi^{v_{j n_j}}, \quad j = 1,\ldots,r.
\]
We set $z_j :=a_{j1} \xi^{v_{j1}} = \cdots = a_{j n_j}\xi^{v_{j n_j}}$. 
It follows from  $v_{j1} + \cdots + v_{j n_j} = 0$ that 
\[
  z_j^{n_j} = a_{j1}\cdots a_{j n_j},\quad j = 1,\ldots,r.
\]
Denote the product $a_{j1}\cdots a_{j n_j}$ by $y_j$. 
It is easy to compute the values of the polynomials $f$ and  $H_f^0(t)$  
at a critical point $\xi$: 
\[
  f(\xi) = 1 - n_1 z_1 - \cdots - n_r z_r,
\]
\[
  H_f^0(\xi) = (-1)^d n_1\cdots n_r z_1^{d_1} \cdots z_r^{d_r}.
\]
On the other hand, we have 
\[P(a_{11}\xi^{v_{11}}, \ldots,a_{rn_r}\xi^{v_{rn_r}} ) = 
  z_1^{k_1}\cdots z_r^{k_r}.
\]
Thus, we obtain
\begin{eqnarray*}
R_{x^k}(u) = \frac{1}{n_1\cdots n_r}
                       \sum_{\begin{subarray}{l}
                              z_1^{n_1} = y_1 \\
                              \cdots          \\
                              z_r^{n_r} = y_r
                             \end{subarray}} 
   \frac{z_1^{k_1}\cdots z_r^{k_r}}
        {(1 - n_1 z_1 - \cdots - n_r z_r)
         z_1^{d_1}\cdots z_r^{d_r}}.
\end{eqnarray*}
Setting  $u_j := n_j^{n_j} y_j$ ($1\le j \le r$), we get the required
formula. 
\end{proof}

In order to compute the generating functions of intersection
numbers $I(P, \beta)$  on the Morrison-Plesser moduli spaces associated 
with $\P$,  we remark that  the
Mori cone $K_{\rm eff}(\P)$ is a simplicial $r$-dimensional cone 
generated by the canonical basis of $H_2(\P,\Z) \cong \Z^r$. 
Moreover,  all Morrison-Plesser moduli spaces $\P_\beta$ 
are projective spaces 
$\P^{n_1 b_1 + d_1}\times \cdots \times \P^{n_r b_r + d_r}$, 
where $\beta$ can be identified with 
a lattice point $(b_1, \ldots, b_r) \in  K_{\rm eff}(\P) = 
\R_{\geq 0}^r$. Therefore, the generating function 
 \[
  I_{P}(y) = \sum_{\beta \in K_{\rm eff}(\P)}  
  I(P,\beta)\, a^{\beta},
\]
can be rewritten as 
\[
  I_{P}(y) = \sum_{b_1,\ldots,b_r\ge 0}  
  I(P,\beta)\, y_1^{b_1}\cdots y_r^{b_r},
\]
where $y_j  = a_{j1}\cdots a_{j n_j}$ $( 1 \leq j \leq r)$. 

\begin{theo} 
\label{gen_prod}
Let $u_j = n_j^{n_j} y_j$ $(1 \le j \le r)$. 
Then the generating function of intersection numbers $I(x^k, \beta)$ 
can be written as 
\[
  I_{x^k}(u) = n_1^{d_1 - k_1}\cdots n_r^{d_r - k_r}
  \sum_{b_1,\ldots,b_r\ge 0} 
  \frac{(n_1 b_1 + \cdots + n_r b_r)!}
       {(n_1 b_1 + d_1 - k_1)!\cdots (n_r b_r + d_r - k_r)!}
  u_1^{b_1}\cdots u_r^{b_r}.
\]
\end{theo}

\begin{proof}
Let $[H_j]$ be the hyperplane class on $\P^{d_j}$. Since $\P$ contains 
exactly $n_j$ torus-invariant divisors having the class $[H_j]$, 
the intersection 
number  $I(x^k,\beta)$ equals 
\[
   \<[H_1]^{k_1}\cdots [H_r]^{k_r} 
  (n_1 [H_1] + \cdots + n_r [H_r])^{n_1 b_1 + \cdots + n_r b_r}\>_\beta.
\]
We have 
\[
  I(x^k,\beta) = 
  \sum_{m_1 + \cdots + m_r = m}
  \frac{m!}{m_1!\cdots m_r!}\, 
  n_1^{m_1}\cdots n_r^{m_r}\, 
  \<[H_1]^{m_1 + k_1}\cdots [H_r]^{m_r + k_r}\>_\beta, 
\]
where $m = \sum_{j = 1}^r n_j b_j$. 
The intersection theory on $\P_\beta$ implies 
\[
  \<[H_1]^{l_1}\cdots [H_r]^{l_r}\>_\beta = 
  \left\{
  \begin{array}{ll}
  1,   & l_j = n_j b_j + d_j, \quad j = 1,\ldots,r, \\
  0,   & {\rm otherwise},
  \end{array}
  \right.
\]
So we obtain 
\[
  I_{x^k}(u) = n_1^{d_1 - k_1}\cdots n_r^{d_r - k_r}
  \sum_{b_1,\ldots,b_r\ge 0} 
  \frac{(n_1 b_1 + \cdots + n_r b_r)!}
       {(n_1 b_1 + d_1 - k_1)!\cdots (n_r b_r + d_r - k_r)!}
  u_1^{b_1}\cdots u_r^{b_r},
\]
where $u_j = n_j^{n_j} y_j$ ($1 \le j \le r$).
\end{proof}
\medskip

It is sufficient  to verify the Toric Residue Mirror Conjecture 
(Conjecture~\ref{TRMC}) for any monomial $P= x^k$ of degree 
$d$. By Theorem \ref{trace_prod_proj} and Theorem \ref{gen_prod}, the equality
$R_{x^k}(u) = I_{x^k}(u)$ follows from 

\begin{prop} 
Let 
\[
  T_k(u) := {\rm tr}_{L/K} \left( \frac{1}{z_1^{d_1 - k_1}
  \cdots z_r^{d_r - k_r}(1 - z_1 - \cdots - z_r)} \right), 
\]
then $T_k$ has the following power series expansion
\[ 
T_k(u)=  n_1\cdots n_r \sum_{b_1,\ldots,b_r\ge 0} 
  \frac{(n_1 b_1 + \cdots + n_r b_r)!}
       {(n_1 b_1 + d_1 - k_1)!\cdots (n_r b_r + d_r - k_r)!}
  u_1^{b_1}\cdots u_r^{b_r},
\] 

\label{pp_MP}
\end{prop}

\begin{proof}
We have the chain of equalities: 
\begin{eqnarray*}
  && T_k(u) =
  \sum_{\begin{subarray}{l}
               z_1^{n_1} = u_1 \\
               \cdots             \\
               z_r^{n_r} = u_r
           \end{subarray}} 
  \frac{1}{z_1^{d_1 - k_1}\cdots z_r^{d_r - k_r}} 
  \sum_{b_1,\ldots,b_r \ge 0} 
  \frac{(b_1 + \cdots + b_r)!}{b_1!\cdots b_r!}
  z_1^{b_1}\cdots z_r^{b_r} \\
&&=
  \sum_{\begin{subarray}{l}
               z_1^{n_1} = u_1 \\
               \cdots             \\
               z_r^{n_r} = u_r
           \end{subarray}} 
  \frac{1}{z_1^{d_1 - k_1}\cdots z_r^{d_r - k_r}} 
  \sum_{b} 
  \frac{(b_1 + \cdots + b_r - \sum_{i = 1}^r (d_i - k_i))!}
       {b_1!\cdots b_r!}
  z_1^{b_1}\cdots z_r^{b_r} \\
&&=
  \sum_{\begin{subarray}{l}
               z_1^{n_1} = u_1 \\
               \cdots             \\
               z_r^{n_r} = u_r
           \end{subarray}} 
  \frac{1}{z_1^{d_1 - k_1}\cdots z_r^{d_r - k_r}} 
  \sum_{b} 
  \frac{(b_1 + \cdots + b_r)!}
       {(b_1 + d_1 - k_1)!\cdots (b_r + d_r - k_r)!}
  z_1^{b_1 + d_1 - k_1}\cdots z_r^{b_r + d_r - k_r} \\
&&=
  \sum_{\begin{subarray}{l}
               z_1^{n_1} = u_1 \\
               \cdots             \\
               z_r^{n_r} = u_r
           \end{subarray}}  
  \sum_{b} 
  \frac{(b_1 + \cdots + b_r)!}
       {(b_1 + d_1 - k_1)!\cdots (b_r + d_r - k_r)!}
  z_1^{b_1}\cdots z_r^{b_r} \\
&&=
  n_1\cdots n_r \sum_{b_1,\ldots,b_r\ge 0} 
  \frac{(n_1 b_1 + \cdots + n_r b_r)!}
       {(n_1 b_1 + d_1 - k_1)!\cdots (n_r b_r + d_r - k_r)!}
  u_1^{b_1}\cdots u_r^{b_r},
\end{eqnarray*}
where in the last row we have used the identity 
\[
  \sum_{z_j^{n_j} = u_j} z_j^{b_j} = 
  \left\{
  \begin{array}{ll}
  n_j  u_j^{b_j}, & b_j = k n_j, \\
  0,   & b_j \ne k n_j, \quad k = 1,2,\ldots
  \end{array}
  \right.
\]
for each $j = 1,\ldots,r$.
\end{proof}

\bigskip

\section{Yukawa  $(d-1)$-point functions for Calabi-Yau hypersurfaces} 
\label{Sect_Yukawa_CY}

Let $\Delta\subset M_\R$ be a reflexive polytope of dimension $d$, 
$A = \{ 0, v_1, \ldots, v_n\}$ a finite subset in $\Delta\cap M$ containing
$0$ and all vertices of $\Delta$, $f(t)$ a Laurent polynomial in the variables
$t_1,\ldots,t_d$ with support in $\Delta$ of the form 
\[
  f(t) := 1 - \sum_{i = 1}^n a_i t^{v_i}.
\]
We denote by $Z_f$ a Calabi-Yau hypersurface defined by the equation $f=0$ in
the torus $\T\cong (\C^*)^d$. Let 
\[
  \Omega := {\bf Res}
  \left(
  \frac{1}{f}\, \frac{d t_1}{t_1}\wedge\cdots\wedge\frac{d t_d}{t_d}
  \right),
\]
be the image of the canonical $d$-form 
\[
  \frac{1}{f}\, \frac{d t_1}{t_1}\wedge\cdots\wedge\frac{d t_d}{t_d}
  \in H^d(\T\backslash Z_f)
\]
under the {\it Poincar\'e residue mapping}
\[
  {\bf Res} : H^d(\T\backslash Z_f) \rightarrow H^{d - 1}(Z_f).
\]

\begin{dfn} 
Assign to each nonzero lattice point $v_i \in A$ a  variable $x_i$. Let 
$Q(x_1,\ldots,x_n)\in\Q[x_1,\ldots,x_n]$ be a homogeneous polynomial of degree
$d - 1$. The $Q$-{\it Yukawa  $(d-1)$-point function} is defined by the
formula 
\begin{equation}
\label{def_of_Yuk}
  Y_Q(a_1,\ldots,a_n) := 
  (-1)^{\frac{(d - 1)(d - 2)}{2}}\frac{1}{(2\pi i)^{d - 1}}
  \int_{Z_f} \Omega \wedge 
  \displaystyle
  Q\left(a_1 \frac{\partial}{\partial a_1},\ldots,
         a_n \frac{\partial}{\partial a_n}\right) \Omega,
\end{equation}
where the differential operators 
$a_1\partial/\partial a_1,\ldots,a_n\partial/\partial a_n$ are determined  by
the Gau\ss-Manin connection. 
\end{dfn}

\begin{rem} The sign $(-1)^{\frac{(d - 1)(d - 2)}{2}}$ in the definition of
the Yukawa $(d-1)$-point function is inherited by the variation of the Hodge
structure (see in \cite[Section~8.6.3]{CK}). The $3$-point Yukawa functions
are also called {\it Yukawa couplings}. 
\end{rem}

\begin{exam} 
Consider the mirror family of Calabi-Yau hypersurfaces in the projective space
$\P^d$ defined by the  Laurent polynomial 
\[
  f(t) = 1 - \sum_{i = 1}^n a_i t^{v_i}, \quad n = d + 1,
\]
where $v_1, \ldots, v_d$  form a basis of the lattice $M$ and  
\[
  v_{d+1} = -(v_1 + \cdots + v_d).
\]
If we set $y := a_1\cdots a_n$, then the $Q$-Yukawa $(d-1)$-point function is
equal to 
\begin{equation}
\label{Proj_space}
  Y_Q(y) =   \frac{n\, Q(1,\ldots,1)}{1 - n^n\, y},\quad \deg Q = d - 1, 
\end{equation}
(see e.g. \cite{JN}). In the particular case $d = 4$ and $Q(x) = x^3$ the
formula (\ref{Proj_space}) gives the Yukawa $3$-point function for mirrors of
Calabi-Yau quintic hypersurfaces in $\P^4$ 
\[
  Y_{x^3}(y) = \frac{5}{1 - 5^5y}, 
\]
which is well-known from  \cite{COGP}.
\end{exam}

In order to establish the relation between Yukawa $(d-1)$-point functions and
toric residues, we need the notion of homogeneous coordinate ring  of a toric
variety  \cite{Cox1}. Let $\Si$ be a complete simplicial fan in $N_\R$
defining a projective simplicial toric variety $\P = \P_\Si$, 
$\Si(1):= \{ e_1, \ldots, e_r\}$ the set of generators of $1$-dimensional
cones in $\Si$, and $z_1, \ldots, z_r$ the corresponding homogeneous
coordinates. The polynomial ring 
\[
  S(\P) := \C[z_1, \ldots, z_r] 
\] 
having a natural grading by ${\rm Cl}(\P)$ is called the 
{\it homogeneous coordinate ring} of $\P$. 

Let $H = \sum_{i=1}^r c_i D_i$ be a big and nef divisor on $\P$. One obtains
the convex polytope $\Delta_H \subset M_\R$ as intersection of $r$ half-spaces
$\< m, e_i \> \geq -c_i$ $(1 \leq i \leq r)$. For any  lattice point $m$ in
$k\Delta_H$, one has $\<m,e_i\> + k c_i \geq 0$ $(1 \leq i \leq r)$. Thus we
can define a mapping $S_{\Delta_H} \to S(\P)$ which sends a monomial 
$t_0^k\, t^m \in S_{\Delta_H}^k$ to the monomial 
$\prod_{i=1}^r  z_i^{\<m,e_i\> + k c_i} \in S(\P)$. This mapping yields an
isomorphism of graded rings 
\begin{equation}
\label{isom_of_rings}
  S_{\Delta_H} \cong \bigoplus_{k = 0}^\infty S(\P)_{k \alpha},
\end{equation}
where $\alpha$ is the class of $H$ in ${\rm Cl}(\P)$. For all $k$, this
isomorphism  identifies the subspace $I_{\Delta_H}^k \subset S_{\Delta_H}^k$
with the image of $S(\P)_{k\alpha - \omega_0}$ in $S(\P)_{k\alpha}$ under the
mapping 
\[
  S(\P)_{k\alpha -  \omega_0}\stackrel{\prod_{i=1}^r z_i}{\longrightarrow} 
  S(\P)_{k\alpha},
\]
where $\omega_0 \in {\rm Cl}(\P)$ is the anticanonical class of $\P$. This
bijection allows to compare our notion of toric residue from
Section~\ref{Sect_Residues} with the definition of toric residue given by
D.~Cox in \cite{Cox2}. 

By the isomorphism (\ref{isom_of_rings}), we identify the regular sequence 
$G = (G_0,\ldots,G_d)$ of elements in $S_{\Delta_H}^1$ with its image in
$S(\P)_\alpha$. Thus, by \cite[Theorem~5.1]{Cox2} (see also 
\cite[Thorem~4.8]{Mavlyutov1}, where the theorem of Cox was extended to the
case when $H$ is big and nef), the toric residue mapping (\ref{Tor_res_map})
coincides with the residue mapping 
\[
  {\rm Res}_G : 
  S(\P)_\rho \rightarrow \C, \quad \rho = (d + 1)\alpha - \omega_0 
\]
considered by Cox and induces the canonical isomorphism
\[
  S(\P)_\rho/\<G_0,\ldots,G_d\>_\rho \cong \C.
\]

\begin{rem} 
In particular,  let $\alpha = \omega_0$, i.e., $\Delta_H$ be a reflexive
polytope (we denote it shortly by $\Delta$), $f(t)$ be a generic Laurent
polynomial with support $\Delta$ defining a Calabi-Yau hypersurface 
$Z_f \subset \T$, and $F_0(z),\ldots,F_d(z) \in S(\P)_{\omega_0}$ are the
images of the following regular sequence in  $ S_{\Delta}^1$:
\[
  t_0 f, \;
   t_0 t_1 \partial f/\partial t_1,\ldots,
   t_0 t_d \partial f/\partial t_d.
\]
Then the isomorphism (\ref{Res_f_T}) from Section~\ref{Sect_Residues} coincides with
\begin{equation}
\label{aff_res}
  S(\P)_\rho/\<F_0,\ldots,F_d\>_\rho \cong \C, 
  \quad \rho = d \omega_0.
\end{equation}
\label{reflex}
\end{rem} 

\begin{dfn}[\cite{Mavlyutov1}] Given a subset 
$I = \{e_{i_0},\ldots, e_{i_d}\}\subset\Si(1)$ consisting of $d + 1$ elements
and an integral basis $m_1,\ldots,m_d$ of the lattice $M$, denote by $c_I$ the 
determinant of the $(d + 1)\times (d + 1)$-matrix obtained from the matrix
$(\<m_j,e_{i_k}\>_{1\le j\le d, i_k\in I})$ by adding the first row
$(1,\ldots,1)$. 
\label{c-i}
\end{dfn}

The following statement shows how Yukawa $(d-1)$-point functions can be
computed by means of toric residues $ {\rm Res}_{F_I}$ with respect to some
sequence $F_I$ of polynomials in the  homogeneous coordinate ring of a toric
variety (it is a reformulation of the result in \cite[p.~104]{Mavlyutov1}): 

\begin{theo} 
\label{Theo_Anv}
Let $\Delta$ be a $d$-dimensional reflexive polytope and 
$F_0(z) \in S(\P)_{\omega_0}$ a generic homogeneous polynomial as in
Remark~\ref{reflex}. Choose a subset 
$I = \{e_{i_0},\ldots,e_{i_d}\}\subset \Si(1)$ is such that $c_I \ne 0$ (see
Definition~\ref{c-i}) and define the sequence $F_I$ of homogeneous polynomials
as $F_I = (z_{i_0}\partial F_0/\partial z_{i_0},\ldots, 
           z_{i_d}\partial F_0/\partial z_{i_d})$. Let 
$Q(x_1,\ldots,x_n)\in\Q[x_1,\ldots,x_n]$ be a homogeneous polynomial of degree
$d - 1$. We set 
\[
  q(t) := (-1)^{d - 1}\, t_0^{d - 1}\,Q(a_1 t^{v_1},\ldots,a_n t^{v_n}) 
\]
and denote by $\tilde{q}(z) \in S(\P)_{(d - 1) \omega_0}$ the image of the
polynomial $q(t)$ under the isomorphism (\ref{isom_of_rings}). Then the
$Q$-Yukawa $(d-1)$-point function is equal to the toric residue 
\[
  Y_Q(a) = - c_I\, {\rm Res}_{F_I}
  \left(\tilde{q}(z)\cdot \prod_{i=1}^r z_i \right). 
\]
\end{theo}

There exists another formula for Yukawa $(d-1)$-point function $Y_Q(a)$ which
does not depend on the choice of a subset 
$I = \{e_{i_0},\ldots,e_{i_d}\}\subset \Si(1)$: 

\begin{theo} 
\label{theo_Yukawa}
Let $Q(x_1,\ldots,x_n)\in\Q[x_1,\ldots,x_n]$ be a homogeneous polynomial of
degree $(d - 1)$ and 
\[
  P(x_1,\ldots,x_n) := (x_1 + \cdots + x_n) Q(x_1,\ldots,x_n).
\]
Then the $Q$-Yukawa $(d-1)$-point function is equal to the toric residue 
\[
  Y_Q(a_1,\ldots,a_n) = 
  (-1)^d\, {\rm Res}_f(t_0^d\, P(a_1 t^{v_1},\ldots, a_n t^{v_n})).
\]
\end{theo}

\begin{proof} 
It follows from the definition of toric residue that 
\begin{eqnarray*}
         {\rm Res}_f (t_0^d\, P(a_1 t^{v_1},\ldots,a_n t^{v_n}))
         &=&  
         {\rm Res}_f (t_0^d\, (1 - f(t))\,Q(a_1t^{v_1},\ldots,a_n t^{v_n})) \\
         &=&  
         {\rm Res}_f (t_0^d\, Q(a_1 t^{v_1},\ldots,a_n t^{v_n})).
\end{eqnarray*}
Let $F_0(z),\ldots,F_d(z)\in S(\P)_{\omega_0}$ be the homogeneous polynomials
as in Remark~\ref{reflex}. Then the last residue can be written in homogeneous 
coordinates as 
\[
  {\rm Res}_f (t_0^d\, Q(a_1 t^{v_1},\ldots,a_n t^{v_n})) =
  (-1)^{d-1}{\rm Res}_F\left(\tilde{q}(z)\cdot\prod_{i = 1}^r z_i\right),
\]
where ${\rm Res}_F$ is the Cox's residue with respect to the sequence 
$F = (F_0,\ldots,F_d)$. Compare this residue with the residue from
Theorem~\ref{Theo_Anv}. To do this, note that the reflexivity of $\Delta$
implies that each monomial $t_0\, t^m \in S_\Delta^1$ maps to the monomial 
$\prod_{i = 1}^r z_i^{\<m,e_i\> + 1}\in S(\P)_{\omega_0}$ in
(\ref{isom_of_rings}). Hence, it is easy to see that the sequences of
homogeneous polynomials $F = (F_0,\ldots,F_d)$ and 
$F_I = (z_{i_0}\partial F_0/\partial z_{i_0},\ldots,
         z_{i_d}\partial F_0/\partial z_{i_d})$ are related by the formulae 
\[
  z_{i_k}\partial F_0/\partial z_{i_k} = 
  F_0 + \sum_{j = 1}^d \<m_j, e_{i_k}\> F_j = \sum_{j = 0}^d A_{j k} F_j,
  \quad k = 0,\ldots,d.
\]
Now, we have  
\begin{eqnarray*} 
  {\rm Res}_f (t_0^d\, P(a_1 t^{v_1},\ldots,a_n t^{v_n})) &=& 
  (-1)^{d-1}{\rm Res}_F\left(\tilde{q}(z)\cdot \prod_{i = 1}^r
  z_i\right) \\ &=&
  (-1)^{d-1}{\rm Res}_{F_I}\left(\det(A_{j k})\, 
   \tilde{q}(z)\cdot \prod_{i = 1}^rz_i\right) \\ &=& 
  (-1)^{d-1}c_I\,
  {\rm Res}_{F_I}\left(\tilde{q}(z)\cdot \prod_ {i = 1}^rz_i\right) = 
  (-1)^d Y_Q(a),
\end{eqnarray*}
where the second row follows from the {\it Global Transformation Law} for
toric residue (see Theorem~\ref{GlTL} below) and the third row follows from
the equality $\det(A_{j k}) = c_I$. 
\end{proof}

Next statement is a particular case of \cite[Theorem~0.4]{CCD}.

\begin{theo}
\label{GlTL}
Let $G = (G_0,\ldots,G_d)$ and $H = (H_0,\ldots,H_d)$ be the regular sequences
of elements in $S(\P)_{\omega_0}$. If 
\[
  H_j = \sum_{i = 0}^d A_{ij} G_i,
\]
where $A_{ij}$ are complex numbers, then for each $P\in S(\P)_{d \omega_0}$,
we have 
\[
  {\rm Res}_G(P) = {\rm Res}_H(P\, \det(A_{i j})).
\]
\end{theo}

\medskip

Toric Residue Mirror Conjecture (Conjecture~\ref{TRMC}) implies the following: 

\begin{coro} 
Let $\Delta \subset M_{\R}$ be an arbitrary reflexive $d$-dimensional polytope
and $A$ a finite subset in $\Delta \cap M$ containing $0$ and all vertices of
$\Delta$. Choose any coherent triangulation 
${\mathcal T} =\{\tau_1, \ldots, \tau_k\}$ of $\Delta$ associated with $A$
such that $0$ is a vertex of all the simplices 
$\tau_1, \ldots, \tau_k$. Denote by $\P = \P_{\Si({\mathcal T})}$ the
simplicial toric variety defined by the fan 
$\Sigma = \Sigma({\mathcal T}) \subset M_\R$ whose $d$-dimensional cones are
exactly $\sigma_i := \R_{\geq 0} \tau_i$ $( 1 \leq i \leq k)$. If 
$A = \{0, v_1, \ldots, v_n\}$ and 
\[ 
  f(t) := 1 - \sum_{i = 1}^n a_i t^{v_i},
\]
then for any homogeneous polynomial $Q(x_1,\ldots,x_n)\in\Q[x_1,\ldots,x_n]$
of degree $d - 1$ the Laurent expansion of the $Q$-Yukawa $(d-1)$-point
function $Y_Q(a)$ at the vertex $v_{\mathcal T} \in {\rm Sec}(A)$
corresponding to the coherent triangulations ${\mathcal T}$ coincides with the
generating function of intersection numbers 
\[
  \sum_{\beta\in K_{\rm eff}(\P)}
  \<Q([D_1],\ldots,[D_n])
  ([D_1] + \cdots + [D_n])\,\Phi_\beta\>_\beta\;a^{\beta},
\]
where the sum runs over all integral points $\beta = (b_1,\ldots,b_n)$ of the
Mori cone $K_{\rm eff}(\P)$, and $a^\beta := a_1^{b_1}\cdots a_n^{b_n}$. 
\end{coro}

\bigskip

\section{Algorithmic calculation of toric residues}

An effective procedure for computing of toric residues in homogeneous
coordinates using Gr\"obner basis calculus was developed in 
\cite{CCD, CD}. Next we describe how these ideas can be used in concrete
calculations of Yukawa $(d-1)$-point functions. We hold the same notations as
in Section~\ref{Sect_Yukawa_CY}.

\begin{dfn}[\cite{Cox1}] 
\label{hom_Jac}
Pick a subset $I = \{i_1,\ldots,i_d\}\subset \{1,\ldots,n\}$ such that 
$e_{i_1},\ldots,e_{i_d}$ are linearly independent. Then define a 
{\it toric Jacobian} for any $(d + 1)$ homogeneous polynomials 
$G_0,\ldots,G_d\in S(\P)_\alpha$ by the formula 
\[
  J_G = 
  \begin{pmatrix} 
  G_0                                   &  
  \cdots                                & 
  G_d                                   \\
  \partial G_0/\partial z_{i_1} & 
  \cdots                                & 
  \partial G_d/\partial z_{i_d} \\
  \vdots &  
  \ddots & 
  \vdots   \\
  \partial G_0/\partial z_{i_d} & 
  \cdots & 
  \partial G_d/\partial z_{i_d}
  \end{pmatrix} / \det(e_I) \widehat{z}_I,
\]
where $e_I = \det (\<m_j,e_{i_k}\>_{1\le j,k\le d})$ and 
$\widehat{z}_I = \prod_{i\not\in I} z_i$.
\end{dfn}

\begin{rem}
Suppose that $G = (G_0,\ldots,G_d)$ is the image in $S(\P)_\alpha$ of the
regular sequence of elements in $S_{\Delta_H}^1$. Then it is easy to show (see
\cite{CDS}) that the Jacobian $J_G$ in Definition \ref{hom_Jac} can be
identified with the image of the Jacobian $J_G$ in (\ref{toric_Jac}) under the
isomorphism (\ref{isom_of_rings}). 
\end{rem}

As we have seen in the proof of Theorem~\ref{theo_Yukawa}, the function 
$Y_Q(a)$ is equal to the Cox's residue 
$-{\rm Res}_F(\tilde{q}(z)\cdot\prod_{i=1}^rz_i)$. We can  compute this
residue using the following method. Choose a Gr\"obner basis of the ideal
generated by $F_0(z),\ldots,F_d(z)\in S(\P)_{\omega_0}$. Then compute the
normal form $normalf(H)$ for the polynomial 
$H = -\tilde{q}(z)\cdot\prod_{i=1}^r z_i$ and the normal form $normalf(J_F)$
for the toric Jacobian $J_F$. Since the quotient in (\ref{aff_res}) is
one-dimensional, both normal forms are the multiples 
\[
  normalf(H) = c\,z^\sigma,\quad normalf(J_F) = c_F\,z^\sigma
\]
of some monomial $z^\sigma\in S(\P)_{d\omega_0}$ by constants $c$ and 
$c_F\ne 0$ modulo $\<F_0,\ldots,F_d\>$. Recall that the toric residue 
${\rm Res}_F(J_F)$ is equal to the normalized volume ${\rm Vol}(\Delta)$ of 
polytope $\Delta$. Now given a polynomial $f(t)$ supported in the reflexive 
polytope $\Delta$ and a homogeneous polynomial $Q$, we get in output 
\[
  Y_Q(a) = \frac{c}{c_F}\,{\rm Vol}(\Delta).
\]

\begin{exam}
We illustrate this method by the following example 
(cf. \cite[Section 8.1]{COFKM}, \cite[Section 4.4]{MP}, 
\cite[Appendix A.1]{HKTY}, \cite[Example 5.6.2.1]{CK}). 
Let $\P(1,1,2,2,2)$ be the weighted projective space defined by the fan 
$\Si'\subset M_\R\cong \R^4$ with one-dimensional generators 
\begin{align*}
  v_1 &= (-1,-2,-2,-2),&  v_2 &= (1,0,0,0), &v_3 = (0,1,0,0), \\ 
  v_4 &= (0,0,1,0),&      v_5 &= (0,0,0,1), &
\end{align*}
which are the vertices of the reflexive polytope $\Delta\subset M_\R$. The
polytope $\Delta$ contains only two lattice points except of listed above: 
$v_6 = (0,-1,-1,-1) = \frac{1}{2}(v_1 + v_2)$ and the origin 
$v_0 = (0,0,0,0)$. Including the additional lattice point $v_6$ to the
generators of $\Si'$ corresponds to the blowup of $\P(1,1,2,2,2)$. The dual
polytope $\Delta^*\subset N_\R$ to $\Delta$ is spanned by the vectors 
\begin{align*}
  e_1 &= (-1,-1,-1,-1),&  e_2 &= (7,-1,-1,-1), &e_3 = (-1,3,-1,-1), \\ 
  e_4 &= (-1,-1,3,-1),&   e_5 &= (-1,-1,-1,3), &
\end{align*}
generating the fan $\Si\subset N_\R$ of $\P = \P_\Si$.

Take the Laurent polynomial 
\[
  f(t) = 1 - \sum_{i = 1}^6 a_i t^{v_i} = 
         1 - a_1 t_1^{-1} t_2^{-2} t_3^{-2} t_4^{-2} - 
             a_2 t_1 - a_3 t_2 - a_4 t_3 - a_5 t_4 - 
             a_6 t_2^{-1}t_3^{-1}t_4^{-1}
\]
having $\Delta$ as support polytope. After choosing the new variables 
\[ 
  y_1 := a_3 a_4 a_5 a_6, \quad 
  y_2 := \frac{a_1 a_2}{a_6^2},
\]
we can put the equation for $f(t)$ in the form 
\[
  f(t) = 1 - y_2 t_1^{-1} t_2^{-2} t_3^{-2} t_4^{-2} - t_1 - y_1 t_2 - 
  t_3 - t_4 - t_2^{-1}t_3^{-1}t_4^{-1}. 
\]
Let $Z_f$ be a hypersurface in $\T \cong (\C^*)^4$ defined by $f(t)$ and
forming the mirror family of Calabi-Yau hypersurfaces in $\P$. Denote 
\[
  \Omega := {\bf Res}
  \left(
  \frac{1}{f}\,\frac{d t_1}{t_1}\wedge\cdots\wedge\frac{d t_4}{t_4}
  \right)
  \in H^3(Z_f).
\]
We compute the Yukawa couplings 
\[
  Y^{(3 - k,k)}(y_1,y_2) = 
  \frac{-1}{(2\pi i)^3}\int_{Z_f} \Omega\wedge 
  (y_1 \partial/\partial y_1)^{3 - k} 
  (y_2 \partial/\partial y_2)^k \Omega,\quad k = 0,1,2,3
\] 
corresponding to polynomials $Q(x_1,x_2) = x_1^{3 - k}x_2^k$. Show how to
compute, say, the Yukawa coupling $Y^{(3,0)}(y_1,y_2)$ corresponding to the
polynomial $Q(x_1,x_2) = x_1^3$. Let $z_i$ be the homogeneous coordinates on
$\P$ related with the vectors $e_i$ ($1\le i \le 5$). Then the
homogenization of $f(t)$ defining the anticanonical hypersurface on $\P$ is 
\[
  F(z) = z_1 z_2 z_3 z_4 z_5  - y_2 z_1^8 - z_2^8 - y_1 z_3^4 - z_4^4 - z_5^4
  - z_1^4 z_2^4.
\]
Denote by $F_1(z),\ldots,F_4(z)$ the images of 
$t_0 t_1\partial f/\partial t_1,\ldots,
 t_0 t_4\partial f/\partial t_4$
under the isomorphism (\ref{isom_of_rings}). Fix a Gr\"obner basis of the
ideal $\<F_0,\ldots,F_4\>$ with respect to the reverse lex order. Next, note
that the homogenization of 
\[
  q(t) = -t_0^3\, Q(y_1 t^{v_1}, y_2 t^{v_3}) = 
  -t_0^3\, Q(y_1 t_2, y_2 t_1^{-1} t_2^{-2} t_3^{-2} t_4^{-2}) = 
  -t_0^3\, (y_1 t_2)^3
\]
is $\tilde{q}(z) = - y_1^3 z_3^{12}$. Applying the Gr\"obner basis
calculation, we have found (using MAPLE) the following normal forms: 
\begin{eqnarray*}
  && normalf(-\tilde{q}(z)\cdot z_1\cdots z_5) = 4y_1^4 z_3^{16}, \\
  && normalf(J_F) = 4y_1^4 ((1 - 2^8 y_1)^2 - 2^{18} y_1^2y_2) z_3^{16}, 
\end{eqnarray*}
where $J_F$ is the toric Jacobian. This easily yields that the Yukawa coupling
is given by the formula 
\[
  Y^{(3,0)}(y_1,y_2) = \frac{8}{(1 - 2^8 y_1)^2 - 2^{18} y_1^2y_2},
\]
since ${\rm Vol}(\Delta) = 8$. Using the same procedure, we obtain 
\begin{eqnarray*}
  Y^{(2,1)}(y_1,y_2) &=& 
  \frac{4(1 - 2^8 y_1)}
       {(1 - 2^8 y_1)^2 - 2^{18} y_1^2y_2}, \\
  Y^{(1,2)}(y_1,y_2) &=& 
  \frac{8y_2(-1 + 2^9 y_1)}
       {(1 - 4 y_2)((1 - 2^8 y_1)^2 - 2^{18} y_1^2y_2)}, 
\\
  Y^{(0,3)}(y_1,y_2) &=&
  \frac{4y_2(1 - 2^8 y_1 + 2^2 y_2 - 2^{10}3 y_1 y_2)}
       {(1 - 4y_2)^2((1 - 2^8 y_1)^2 - 2^{18} y_1^2y_2)}.
\end{eqnarray*}

The series expansions for $ Y^{(3 - k,k)}(y_1,y_2)$ were computed in 
\cite[Section 4.4]{MP}. The result of these computations is 
\[
  Y^{(3 - k,k)}(y_1,y_2) = 
  \sum_{\lambda_1,\lambda_2 \ge 0}
  2^{8\lambda_1 + 2\lambda_2 + 3 - k}
  \binom{\lambda_1 + 1 - k}{2\lambda_2 + 1 - k} \, 
  y_1^{\lambda_1}y_2^{\lambda_2}, \quad k = 0,1,2,3. 
\]

Now suppose that we do not blowup $\P(1,1,2,2,2)$. It means that we do not add
$v_6$ to the generators of $\Si'$. Then the polynomial 
\[
  f(t) = 1 - \sum_{i = 1}^5 a_i t^{v_i} =  
  1 - a_1 t_1^{-1} t_2^{-2} t_3^{-2} t_4^{-2} - a_2 t_1 - a_3 t_2 - 
  a_4 t_3 - a_5 t_4
\]
can be transformed to 
\[
  f(t) = 1 - t_1^{-1} t_2^{-2} t_3^{-2} t_4^{-2} - y t_1 - t_2 - 
  t_3 - t_4
\]
after introducing the new variable  $y := a_1 a_2 a_3^2 a_4^2 a_5^2$. It is
nice to observe that the Yukawa coupling 
\[
  Y^3(y) = \frac{8}{1 - 2^{18} y} = 
  8 + 2097152 y + 549755813888 y^2 + O(y^3)
\]
for $Q(x) = x^3$ can be computed either by theorems~\ref{wp_MP} and 
\ref{theo_Yukawa} or as a limit $a_6 \rightarrow 0$ in $Y^{(3,0)}(y_1,y_2)$ 
found above. 

\end{exam}

\bigskip

\section{Mirrors of Calabi-Yau hypersurfaces in $\P^{d_1}\times \P^{d_2}$} 

Let us illustrate our method for computation of  Yukawa $(d-1)$-point
functions for Calabi-Yau hypersurfaces $Z_f \subset \T\cong (\C^*)^2$ defined
by polynomial 
\[
  f(t) = 1 - \sum_{i_1 = 1}^{n_1} a_{1 i_1} t^{v_{1 i_1}} - 
             \sum_{i_2 = 1}^{n_2} a_{2 i_2} t^{v_{2 i_2}},
\]
where the exponents $v_{ij}$ satisfy the relations
\begin{eqnarray*}
  v_{11} + \cdots + v_{1n_1} = 0, \quad
  v_{21} + \cdots + v_{2n_2} = 0.
\end{eqnarray*}
The toric variety $\P_\Si$ with the fan $\Si$ of dimension $d = d_1 + d_2$
(here, $d_1 = n_1 - 1, d_2 = n_2 - 1$) with generators $\{v_{j n_j}\}$ is the
product of two projective spaces $\P^{d_1}\times\P^{d_2}$. According to
\cite{Batyrev2}, the  hypersurfaces $Z_f$ form  the mirror family of
Calabi-Yau hypersurfaces in $\P^{d_1}\times \P^{d_2}$. 

Denote by $y_1 = a_{11}\cdots a_{1n_1}$, $y_2 = a_{21}\cdots a_{2n_2}$ the new
variables and by $\theta_1 = y_1\partial/\partial y_1$, 
$\theta_2 = y_2\partial/\partial y_2$ the logarithmic partial
derivations. Given a form-residue 
\[
  \Omega := {\bf Res}
  \left(
  \frac{1}{f}\,\frac{dt_1}{t_1}\wedge\cdots\wedge\frac{dt_d}{t_d}
  \right)
  \in H^{d - 1}(Z_f), 
\]
the $2$-parameter Yukawa $(d-1)$-point functions are defined as the integrals 
\begin{equation*}
  Y^{(k_1,k_2)}(y_1,y_2) = 
  \frac{(-1)^{\frac{(d - 1)(d - 2)}{2}}}{(2\pi i)^{d - 1}}
  \int_{Z_f}\Omega\wedge \theta_1^{k_1} \theta_2^{k_2}\Omega,
  \quad k_1 + k_2 = d - 1.
\end{equation*}

Theorem~\ref{theo_Yukawa}, Theorem~\ref{trace_prod_proj} and
Proposition~\ref{pp_MP} easily implies the following statement. 

\begin{prop} The Yukawa $(d-1)$-point function $Y^{(k_1,k_2)}(y_1,y_2)$ is 
equal to the toric residue related with the polynomial 
$P(x_1,x_2) = x_1^{k_1} x_2^{k_2}(n_1 x_1 + n_2 x_2)$. Let 
$u_1 = n_1^{n_1} y_1$,  $u_2 = n_2^{n_2} y_2$. Then the function 
$Y^{(k_1,k_2)}(y_1,y_2)$ can be computed 

\begin{enumerate}
\item[1)] as a trace of rational function: 
\begin{eqnarray*}
  Y^{(k_1,k_2)}(u_1, u_2) = 
  n_1^{d_1 - k_1 - 1} n_2^{d_2 - k_2 - 1}
  \sum_{\begin{subarray}{l}
                    z_1^{n_1} = u_1 \\ 
                    z_2^{n_2} = u_2 
                   \end{subarray}}
  \frac{1}{z_1^{d_1 - k_1} z_2^{d_2 - k_2}(1 - z_1 - z_2)};
\end{eqnarray*}

\item[2)] as a series: 
\[
  Y^{(k_1,k_2)}(u_1,u_2) = 
  n_1^{d_1 - k_1} n_2^{d_2 - k_2}\sum_{b_1,b_2\ge 0} 
  \frac{(n_1 b_1 + n_2 b_2 + 1)!}
  {(n_1 b_1 + d_1 - k_1)!(n_2 b_2 + d_2 - k_2)!} u_1^{b_1} u_2^{b_2}.
\]
\end{enumerate}
\end{prop}

\medskip

Some examples of explicit calculation of Yukawa $(d-1)$-point functions as
rational functions are given below. 

\begin{exam}
$\P^1\times\P^1$; $u_1 = 2^2 y_1$, $u_2 = 2^2 y_2$.
\begin{eqnarray*}
  Y^{(1,0)}(u_1,u_2) &=& \frac{1}{2}
  \sum_{\begin{subarray}{l}  
          z_1^2 = u_1 \\
          z_2^2 = u_2
        \end{subarray}}
  \frac{1}{z_2(1 - z_1 - z_2)} = 
  2 \sum_{b_1,b_2\ge 0}
  \frac{(2b_1 + 2b_2 + 1)!}{(2b_1)!(2b_2 + 1)!} u_1^{b_1} u_2^{b_2} \\
  &=&
  \frac{2(1 + u_1 - u_2)}{(1 - u_1 - u_2)^2 - 4 u_1 u_2}.
\end{eqnarray*}
By symmetry, $Y^{(0,1)}(u_1,u_2) = Y^{(1,0)}(u_2,u_1)$.
\end{exam}

\begin{exam}
$\P^1\times\P^2$; $u_1 = 2^2 y_1$, $u_2 = 3^3 y_2$.
\begin{eqnarray*}
  Y^{(2,0)}(u_1,u_2) &=& \frac{3}{2}
  \sum_{\begin{subarray}{l}  
          z_1^2 = u_1 \\
          z_2^3 = u_2
        \end{subarray}}
  \frac{z_1}{z_2^2(1 - z_1 - z_2)} = 
  \frac{9}{2} \sum_{b_1,b_2\ge 0}
  \frac{(2b_1 + 3b_2 + 1)!}{(2b_1 - 1)!(3b_2 + 2)!} u_1^{b_1} u_2^{b_2} \\
  &=&
  \frac{9}{2}
  \frac{u_1(3 + u_1)}{(1 - u_1)^3 - 2 u_2 (1 - 3 u_1)}.
\end{eqnarray*}

\begin{eqnarray*}
  Y^{(1,1)}(u_1,u_2) &=& \frac{1}{2}
  \sum_{\begin{subarray}{l}  
          z_1^2 = u_1 \\
          z_2^3 = u_2
        \end{subarray}}
  \frac{1}{z_2(1 - z_1 - z_2)} = 
  3 \sum_{b_1,b_2\ge 0}
  \frac{(2b_1 + 3b_2 + 1)!}{(2b_1)!(3b_2 + 1)!} u_1^{b_1} u_2^{b_2} \\
  &=&
  \frac{3(1 - u_2 - u_1^2)}{(1 - u_1)^3 - 2 u_2 (1 - 3 u_1)}.
\end{eqnarray*}

\begin{eqnarray*}
  Y^{(0,2)}(u_1,u_2) &=& \frac{1}{3}
  \sum_{\begin{subarray}{l}  
          z_1^2 = u_1 \\
          z_2^3 = u_2
        \end{subarray}}
  \frac{1}{z_1(1 - z_1 - z_2)} = 
  2 \sum_{b_1,b_2\ge 0}
  \frac{(2b_1 + 3b_2 + 1)!}{(2b_1 + 1)!(3b_2)!} u_1^{b_1} u_2^{b_2} \\
  &=&
  \frac{2((1 - u_1)^2 + 2 u_2)}{(1 - u_1)^3 - 2 u_2 (1 - 3 u_1)}.
\end{eqnarray*}
\end{exam}

The next example of hypersurfaces in product of projective spaces was
considered in \cite{BvS} in the case of diagonal one-parameter subfamily and
in \cite{HKTY}. 

\begin{exam}
$\P^2\times\P^2$; $u_1 = 3^3 y_1$, $u_2 = 3^3 y_2$.
\begin{eqnarray*}
  Y^{(3,0)}(u_1,u_2) &=& \frac{1}{3}
  \sum_{\begin{subarray}{l}  
          z_1^3 = u_1 \\
          z_2^3 = u_2
        \end{subarray}}
  \frac{z_1}{z_2^2(1 - z_1 - z_2)} = 
  3 \sum_{b_1,b_2\ge 0}
  \frac{(3b_1 + 3b_2 + 1)!}{(3b_1 - 1)!(3b_2 + 2)!} u_1^{b_1} u_2^{b_2} \\
  &=&
  \frac{9 u_1(2 + u_1 + u_2)}{(1 - u_1 - u_2)^3 - 27 u_1 u_2}.
\end{eqnarray*}

\begin{eqnarray*}
  Y^{(2,1)}(u_1,u_2) &=& \frac{1}{3}
  \sum_{\begin{subarray}{l}  
          z_1^3 = u_1 \\
          z_2^3 = u_2
        \end{subarray}}
  \frac{1}{z_2(1 - z_1 - z_2)} = 
  3 \sum_{b_1,b_2\ge 0}
  \frac{(3b_1 + 3b_2 + 1)!}{(3b_1)!(3b_2 + 1)!} u_1^{b_1} u_2^{b_2} \\
  &=&
  \frac{3((1 - u_2)^2 + u_1(1 - 2 u_1 - u_2))}{(1 - u_1 - u_2)^3 - 27 u_1 u_2}.
\end{eqnarray*}
By symmetry, $Y^{(k_1, k_2)}(u_1,u_2) = Y^{(k_2, k_1)}(u_2,u_1)$.
\end{exam}


\end{document}